\documentclass[11pt]{article}
\usepackage[margin=1in]{geometry}
\usepackage{amsmath,amssymb,amsthm}
\usepackage{graphicx}
\usepackage{array}
\usepackage{xcolor}
\usepackage[colorlinks=true,linkcolor=blue,citecolor=blue,urlcolor=blue]{hyperref}

\theoremstyle{remark}
\newtheorem{remark}{Remark}

\newcommand{\Tr}{\operatorname{Tr}}

\allowdisplaybreaks

\title{{Optimal isotropic three-material conducting composites with an ideal
conductor, in the plane and in space.\\ Bounds and optimal structures}}
\author{Andrej Cherkaev\thanks{Department of Mathematics, University of Utah} 
and Leonid Gibiansky\thanks{QuantPharm LLC)}}
\date{November 2026}

\begin{document}

\maketitle

\begin{abstract}
{The paper describes an optimal isotropic three-dimensional, three-material conducting composite. One material is an ideal conductor, the volume fractions are prescribed, and the objective is to minimize the effective conductivity over all admissible layouts.

We first derive a lower bound for the effective conductivity as a function of the material properties and volume fractions. The composite energy is bounded below by a constrained polyconvex envelope, where quadratic and cubic null Lagrangians are imposed through Lagrange multipliers. We further assume that the fields in optimal structures are positive semidefinite. Together with the use of the cubic null Lagrangian, this assumption makes the bound exact. The resulting bound determines the extremal conductivity and the set of associated support fields in optimal structures. Depending on the parameters, it has five regions with different analytic expressions.

We then prove the bound attainability by constructing optimal structures that realize the bound in every region. They are spherical atoms: a core of the intermediate-conducting material, surrounded by a mantle of wrapped optimal planar structures. In the two extreme regions, a crust is added from the material whose prescribed volume fraction is excessive.}
\end{abstract}
\newpage
\tableofcontents
\newpage


\section{Introduction}\label{sec:intro}

The question of extremal effective conductivity of a mixture is as old as homogenization. For
two isotropic conductors in prescribed proportions, the arithmetic and harmonic means [Wiener 1912]
bound the effective tensor in every direction, and laminates loaded along and across the layers
attain them. For isotropic mixtures, Hashin and Shtrikman [1962] obtained much sharper bounds
by a variational principle with a homogeneous comparison medium. The assemblage
of coated spheres attains them [Hashin 1962; Hashin and Shtrikman 1962].

The two-material problem was solved in the 1980s. Lurie and Cherkaev [1984,
1986] and, independently, Murat and Tartar [1985] (see also Tartar [1985]) found the $Gm$-closure,
the set of all isotropic and anisotropic effective tensors attainable with given volume fractions. Their tool, the
translation method, goes back to Tartar [1979] and Lurie and Cherkaev [1982]. It adds to the energy
a quasiaffine function of the field (a null Lagrangian), whose average depends only on the average
field, and bounds the translated energy by convexity.

The optimal structures are coated ellipsoids or matrix laminates, of rank two in the plane and rank three
in space. Differential and iterative schemes [Lurie and Cherkaev 1985; Norris 1985; Milton 1986; Cherkaev 2000, \S7.3.1;
Milton 2002, \S9.9; Albin and Cherkaev 2006] reach the same effective tensors by adding
infinitesimal layers to a growing composite. For two materials the theory is complete.

With three and more materials the picture changes. The three-material Hashin--Shtrikman bound still
holds. Milton [1981] showed that in general it is not attainable, and that coated spheres attain it
as long as it does not exceed the conductivity of the intermediate material; their cores are made of
either of the better conductors, and their coats of the poorest one. Lurie and Cherkaev [1985]
attained the bound in the same range with multi-coated spheres. Talbot, Willis and Nesi [1995] and
Albin, Conti and Nesi [2007] obtained improved bounds for three and more materials. Cherkaev and
Gibiansky [1996] constructed, by ``imitation'', planar three-material laminates of rank three with
one principal conductivity at the arithmetic or harmonic mean. Unlike for two materials, the
$Gm$-closure therefore contains segments of the sides of the Wiener box.

The planar problem with an ideal conductor has been solved. Nesi [1995] improved the
Hashin--Shtrikman bound by adding to the translation method the pointwise inequality
$\det\nabla u\ge0$ for the fields, a property of every planar conducting composite [Alessandrini
and Nesi 2001; Bauman, Marini and Nesi 2001]. Gibiansky and Sigmund [2000] found planar
three-material conducting structures that, for an ideal conductor, attain the Hashin--Shtrikman
bound in its whole range of attainability [Milton 2002, eq.~(23.8) and \S25.6]. Cherkaev [2000,
ch.~12] analyzed the problem by necessary conditions on the fields of optimal structures, and Albin,
Cherkaev and Nesi [2007] built multirank laminates for the Hashin--Shtrikman range.

Cherkaev [2009] made this approach systematic for several materials. He added to the translation
method pointwise restrictions on the fields of optimal composites and called the approach
``localized polyconvexity''. For an ideal conductor he found the bound in closed form. The bound has
three regimes with different formulas that define specific optimal fields in each material; in two
of them it coincides with the bound of Nesi. He also showed that laminates attain the bound in every
regime, and represented the fields of a laminate by a path in the plane of field eigenvalues
[Cherkaev 2009, Fig.~6; Cherkaev and Zhang 2011, \S3.3]. Cherkaev and Zhang [2011] extended the
planar analysis to anisotropic three-material composites. Finally, Cherkaev [2012] attained the
isotropic bound with radial structures, the ``wheel'' assemblages of radial spikes.

In three dimensions the planar argument does not carry over. Briane, Milton and Nesi [2004]
constructed periodic composites whose corrector determinant changes sign on a set of positive
measure, so Nesi's inequality fails for general three-dimensional structures. Briane and Nesi [2004]
showed that it holds in every laminate, and used it to show that laminates do not attain the
Hashin--Shtrikman bound for small fractions of the poorest material. They also proved Nesi's bound
in space for ``quasiorthogonal'' composites, those in which the second invariant of the field is
nonnegative in the poorest material, among them the coaxial laminates.

Albin and Cherkaev [2006] constructed ``hairy spheres'' that attain the three-material
Hashin--Shtrikman bound in a range of fractions. A hairy sphere has a core of the intermediate
material, a layer of the poorest material with radial conical hairs of the best conductor, and,
optionally, an outer coat of the poorest material. Liu [2011] derived, with a cubic null Lagrangian,
a necessary condition on the volume fractions, which shows that this range is exact for three
materials. He also gave another explicit construction that attains the bound in this whole range,
coated spheres with a shell of periodic E-inclusions [Liu, James and Leo 2007]. Below this range
the Hashin--Shtrikman bound is not attained, and what replaces it remained open. In this paper we
treat the three-dimensional problem in the case where the third material is an ideal conductor,
$k_3=\infty$.

\paragraph{Plan.} Section~\ref{sec:method} states the problem and describes the methods.
Section~\ref{sec:2d} revisits the planar problem and Section~\ref{sec:3d} solves the three-dimensional one.
Section~\ref{sec:discussion} discusses the results and the open questions, and describes the use
of large language models.

\section{The problem and the method}\label{sec:method}

\subsection{Statement of the problem}\label{sec:problem}

Consider a periodic cell $O\subset\mathbb R^d$, $d=2$ or $d=3$, of unit volume filled with three isotropic materials
$S_1$, $S_2$, $S_3$ of conductivities $k_1<k_2<k_3=\infty$. The
material $S_i$ occupies the subdomain $O_i$ with the volume fraction $m_i$, where
$m_3=1-m_1-m_2$. Denote $\langle f\rangle=\int_O f\,dx$, and call $\chi_i$ the characteristic
function of the subdomain $O_i$, so that the conductivity and the volume fractions are
\begin{equation}
k(x)=\sum_{i=1}^3k_i\chi_i(x),\qquad\langle\chi_i\rangle=m_i .
\label{1.1}
\end{equation}

We apply $d$ orthogonal unit loadings in turn; let $u_j$ solve the periodic conductivity problem
\begin{equation}
\nabla\cdot k(x)\nabla u_j=0,\qquad u_j-x_j\ \text{is $O$-periodic},\qquad j=1,\dots,d .
\label{eq:cell}
\end{equation}
The gradients, collected as columns, form the $d\times d$ field matrix, whose average is the unit matrix,
\begin{equation}
E(x)=\bigl(\nabla u_1|\cdots|\nabla u_d\bigr),\qquad \nabla\times E=0,\qquad \langle E\rangle=I,\qquad
E=0\ \text{in }O_3.
\label{eq:E}
\end{equation}

The energy of the $d$ loadings is $\Tr K_{\rm eff}$, which equals $d\,k_{\rm eff}$ for an isotropic
composite. Here $K_{\rm eff}$ and $k_{\rm eff}$ are the effective conductivities of anisotropic and isotropic materials, respectively.
The optimal composite minimizes this energy,
\begin{equation}
d\,k_0=\inf_{\chi_1,\chi_2}\ \min_{E}\ \bigl\langle(k_1\chi_1+k_2\chi_2)|E|^2\bigr\rangle ,
\label{eq:problem}
\end{equation}
over the fields \eqref{eq:E} and over the layouts \eqref{1.1}; here $|E|^2=\Tr(E^TE)$. We denote
by $k_0$ the minimal effective conductivity. The infimum is taken over all layouts, isotropic or not, and therefore $k_0\le k=\tfrac1d\Tr K_{\rm eff}$ for every composite.
The ideal conductor $S_3$ enters only through the volume it occupies and the connectivity it provides;
hence $k_0$ depends on $k_1$, $k_2$ and the fractions $m_1$, $m_2$ only.

\paragraph{Lower and upper bounds for $k_0$.}
We find a pair of inequalities bounding $k_0$ from below and from above.
The lower bound comes from the geometry-independent inequality $k_0\ge k_*$, valid for all
layouts \eqref{1.1}. The upper bound is the conductivity $k_{\rm str}$ of a particular structure.
If the inequalities
\begin{equation}
k_* \leq k_0 \leq k_{\rm str}
\label{eq:dual}
\end{equation}
become equalities, the optimal conductivity $k_0$ is determined.

\paragraph{Nonconvex variational problem.} The problem can be formulated as a nonconvex variational
problem for $E$. We assign Lagrange multipliers (material costs) $\gamma_i$ to the constraints
$\langle\chi_i\rangle=m_i$, with $\gamma_3=0$ and $\gamma_1>\gamma_2$, and exclude $\chi_i$. The result is
the problem
\begin{equation}
\inf_E\ \langle F(E)\rangle,\qquad \nabla\times E=0,\quad \langle E\rangle=I,\quad E\ \text{is $O$-periodic},
\label{eq:multiwell}
\end{equation}
\[
F=\min_{i=1,2,3}F_i,\qquad F_i(E)=k_i|E|^2+\gamma_i\ \ (i=1,2),\qquad
F_3(E)=\begin{cases}0,& E=0,\\ \infty,& E\ne0.\end{cases}
\]
The Lagrangian $F$ is the minimum of three paraboloids; it is nonconvex in a neighborhood of zero of size $|E|<C_0$. Its minimizers $E$ are ordered and bounded,
\begin{equation}
|E_3|<|E_2|<|E_1|\le C_0 .
\label{eq:constraints}
\end{equation}
This implies that the three materials cannot meet at an isolated line, because the fields in them
cannot be singular.

\subsection{Null Lagrangians and integral restrictions}\label{sec:nl}

The lower bound $k_*\le k_0$ is obtained by the \emph{constrained translation method}, or the
constrained polyconvex envelope of $F$. First, we replace the differential constraint
$\nabla\times E=0$ by its integral corollaries, the null Lagrangians.

{\em A function $f(E)$ is called a null Lagrangian, or quasiaffine, if its average over the
cell depends only on the average of $E$, whatever the geometry. }
For the $d\times d$ gradient matrices
the null Lagrangians are spanned by the minors of every order [Morrey 1966; Reshetnyak 1967; Ball,
Currie and Olver 1981]. For the isotropic problem we use the rotation invariants among them, the sums of the principal $p\times p$ minors. Explicitly
\begin{align}
d=2:&\qquad \langle\Tr E\rangle=2,\qquad \langle\det E\rangle=1, \label{eq:nl2d}\\
d=3:&\qquad \langle\Tr E\rangle=3,\qquad \langle\Phi(E)\rangle=6,\qquad \langle\det E\rangle=1,
\label{eq:nl3d}
\end{align}
where
\[
\Phi(E)=(\Tr E)^2-\Tr(E^2)
\]
is twice the sum of the $2\times2$ principal minors.
Since $E=0$ in $O_3$, every structure satisfies the $d$ \emph{integral restrictions}
\begin{align}
d=2:&\qquad m_1\langle\Tr E\rangle_1+m_2\langle\Tr E\rangle_2=2,\qquad
m_1\langle\det E\rangle_1+m_2\langle\det E\rangle_2=1, \label{eq:restr2d}\\
d=3:&\qquad m_1\langle\Tr E\rangle_1+m_2\langle\Tr E\rangle_2=3,\qquad
m_1\langle\Phi\rangle_1+m_2\langle\Phi\rangle_2=6, \notag\\
&\qquad m_1\langle\det E\rangle_1+m_2\langle\det E\rangle_2=1, \label{eq:restr3d}
\end{align}
where $\langle\cdot\rangle_i$ is the average over $O_i$. These are the restrictions on the trace,
 on the second invariant (in space) and on the determinant.

\subsection{The bound: translation with constraints}\label{sec:translation}

The differential constraint is replaced by the integral restrictions \eqref{eq:restr2d}--\eqref{eq:restr3d}. The
replacement does not decrease the set of admissible fields; thus the solution of the resulting
finite-dimensional problem gives a lower bound $k_*$ of $k_0$,
\[
d\,k_*=\min\ \bigl(m_1\langle k_1|E|^2\rangle_1+m_2\langle k_2|E|^2\rangle_2\bigr),\qquad k_*\le k_0 ,
\]
where the minimum is taken over the fields that satisfy  \eqref{eq:restr2d}--\eqref{eq:restr3d} and the pointwise
restriction below.
We add the restrictions with the Lagrange multipliers $r$ (on the trace), $\tau$ (on $\Phi$) and,
in space, $\mu$ (on the determinant). For any structure,
\begin{equation}
d\,k\ \ge\ V(r,\tau,\mu)=\sum_{i=1,2}m_i\,g_i(r,\tau,\mu)+d\,r-d(d-1)\,\tau-\mu ,
\label{eq:wd}
\end{equation}
\begin{equation}
g_i(r,\tau,\mu)=\min_{E\succeq0}\bigl[k_i|E|^2-r\Tr E+\tau\,\Phi(E)+\mu\det E\bigr],
\label{eq:gi}
\end{equation}
and $d\,k_*=\max_{r,\tau,\mu}V$. Notice that the minimum over $E$ is achieved at symmetric matrices, because the antisymmetric part does not change the trace and only increases the other terms. Indeed, write $E=S+W$, with $S$ symmetric and $W$ skew with the axial vector $w$. Then $|E|^2=|S|^2+|W|^2$,
$\Tr E=\Tr S$, $\Phi(E)=\Phi(S)+|W|^2$ and $\det E=\det S+w\cdot Sw$, so the bracket in
\eqref{eq:gi} exceeds its value at $S$ by $(k_i+\tau)|W|^2+\mu\,w\cdot Sw\ge0$ when $S\succeq0$,
$\mu\ge0$ and $\tau\ge-k_i$. In the plane, where only $\det E\ge0$ is required
(Section~\ref{sec:H}), the same follows from the singular values $\sigma_1,\sigma_2$ of $E$:
$|E|^2=\sigma_1^2+\sigma_2^2$ and $\det E=\sigma_1\sigma_2$ depend only on them, while
$\Tr E\le\sigma_1+\sigma_2$, with equality if and only if $E$ is symmetric and nonnegative, so for
$r>0$ the minimizer is symmetric and nonnegative. Therefore, optimal fields $E_i$ are symmetric, and
their eigenvalues $\lambda_i$ are real and nonnegative.

\subsection{Pointwise restriction \texorpdfstring{$E\succeq 0$}{E >= 0}: Conjecture (H)}\label{sec:H}

Notice that in  \eqref{eq:gi}  we restrict the fields to the cone of positive semidefinite matrices,
\begin{equation}
E(x)\succeq0\qquad\text{for almost every }x\text{ in }O_1\cup O_2 ,
\label{eq:psd}
\end{equation}
where, for a nonsymmetric field, \eqref{eq:psd} refers to its symmetric part. In the plane problem
\eqref{eq:psd} may be replaced by the inequality $\det E\ge0$, which every field satisfies (Section~\ref{sec:intro}). In space $\det E\ge0$ is not a property of arbitrary structures (Section~\ref{sec:intro}), and we use \eqref{eq:psd} as a conjecture, (H), discussed below.

\paragraph{Conjecture (H).}

We conjecture that \eqref{eq:psd} is satisfied in all optimal structures.
This is a property of optimal fields, not of every admissible field.
We do not rigorously prove (H) here. The following two arguments support it;
the first is informal, and the second, for the weaker inequality $\det E\ge0$, is proved for laminates.
\begin{enumerate}
\item The three isotropic points that support the convex envelope of the multiwell energy, $E_3=0$,
$\langle E_2\rangle_2 =c_2I$ and $\langle E_1\rangle_1 =c_1I$, cannot coexist in a connected structure.
The energy $\mathcal E$ of any three-material composite with $E_3=0$ can be presented as
$$
\mathcal E= m_1k_1|\langle E_1\rangle_1|^2 + m_2k_2|\langle E_2\rangle_2|^2 + {\it N}_\eta, \quad {\it N}_\eta=m_1k_1\langle |\eta_1|^2\rangle_1 + m_2k_2\langle |\eta_2|^2 \rangle_2
$$
where $\eta_1$ and $ \eta_2$ are deviations from their average values, $ E_i(x) =\langle E_i\rangle_i + \eta_i(x)$. Notice that $\langle \eta_1\rangle_1= \langle \eta_2\rangle_2=0 $.

The deviations $\eta_1$ and $ \eta_2$ depend on the composite structure, and the optimal structure corresponds to the minimum of $ {\it N}_\eta$.   The norm $ {\it N}_\eta$ is not zero  in a connected structure because of the differential constraint $\nabla \times E =0$, which requires that across an interface the
tangential components of a curl-free field are continuous. In particular, a field adjacent to the ideal conductor is rank one, $E=n\otimes v$, with $n$ the normal to the interface.

The support sets (Section~\ref{sec:support}) must therefore contain rank-one fields
\begin{equation}
E_i :   \langle E_i\rangle_i + \eta_i=n\otimes v, \quad i=1, 2
\label{eq:compatib}
\end{equation}
but compatibility alone does not say which. With $v$ arbitrary, the symmetric part of $n\otimes v$ may be indefinite, and these fields together with the isotropic
points span all matrices.

The restriction \eqref{eq:psd} keeps the fields $\lambda\,n\otimes n$,
$\lambda>0$, which are the fields of the layers between faces of the conductor in every structure
built here. The smallest convex cone that contains these fields and the positive isotropic points $\langle E_i\rangle_i$ is the cone of positive semidefinite matrices. This restriction is sufficient for the compatibility.

\item Laminates never change the sign of $\det E$, whatever the normals of their laminations [Briane
and Nesi 2004, Th.~2.13]. It follows that in any laminate of finite rank made of materials of finite conductivities,
with  the mean field $I$, $\det E>0$ in every layer; with an ideal conductor, $\det E\ge0$.

%

\end{enumerate}

Every structure constructed in this paper satisfies \eqref{eq:psd}.

\subsection{Support sets and  optimal layouts}\label{sec:support}
The energy minimum is attained at symmetric fields (Section~\ref{sec:translation}), so it can be
written in the eigenvalues $\lambda\in\mathbb R^d_+$
of $E$,
\begin{equation}
g_i=\min_{\lambda\ge0}\Bigl[k_i\sum_j\lambda_j^2-r\sum_j\lambda_j+2\tau\sum_{j<l}\lambda_j\lambda_l
+\mu\prod_j\lambda_j\Bigr]
\label{eq:gilambda}
\end{equation}
($\mu=0$ in the plane). The translated density $k|E|^2+\tau\Phi(E)-r\Tr E$ has the same form in the
plane and in space; space adds the term $\mu\det E$.


We call the sets of minimizers of \eqref{eq:gilambda} at the optimal multipliers the \emph{support
sets} $\Lambda_1$, $\Lambda_2$ of the two ordinary materials. We write them as sets of eigenvalue
$d$-tuples; for the conductor, $\Lambda_3=\{0\}$. 
%
%

%

Notice that the bound does not prescribe a structure, only a set of admissible
field values $\Lambda_i$ in each material. Every structure whose fields stay in these sets is
optimal; for the classical translation this is the attainability criterion of Milton [1990] and
Grabovsky [1996].

\subsection{Optimal structures}\label{sec:structures}

The bound $k_*$ restricts $k_0$ from below. The dual, upper bound $k_{\rm str}\ge k_0$ comes from the
effective properties of particular structures. Since $V(r,\tau,\mu)\le d\,k_*$ for all
multipliers, a structure with $d\,k_{\rm str}=V(r,\tau,\mu)$ for some multipliers turns both
inequalities \eqref{eq:dual} into equalities: the structure is optimal, and these multipliers
maximize $V$. This is how each bound below is confirmed by a structure.
A structure attains the bound if and only if its fields lie in the
support sets. In Sections~\ref{sec:2d} and~\ref{sec:3d} we use assemblages of coated spheres (discs).
Let us recall how to compute their effective properties.

\paragraph{Assemblages of coated spheres (discs).} Effective medium theory defines the effective
conductivity as the conductivity of a medium in which an inhomogeneous (composite)
sphere (disc) does not perturb a homogeneous field outside. The inclusion, which we call an \emph{atom}, is
then perfectly hidden in the effective medium, and an assemblage of such atoms of all sizes filling
the space has this conductivity [Hashin 1962; Hashin and Rosen 1964]. 

Consider a spherical inclusion
of unit radius that consists of a core of conductivity $k_c$ surrounded by a coat of conductivity
$k_j$ that takes the fraction $\beta$ of its volume. In a homogeneous field in the medium of
conductivity $k$ it perturbs the potential outside it as
\[
u=\Bigl(1-\frac{\Delta}{\rho^d}\Bigr)\rho\cos\varphi ,
\]
where $\rho$ is the distance from the center, $\varphi$ the angle to the field, and $\Delta$ depends on
$k$, $k_j$, $k_c$ and $\beta$. There is no perturbation, $\Delta=0$, when the parameters are related as 
\begin{equation}
\frac{1}{k+(d-1)k_j}=\frac{\beta}{d\,k_j}+\frac{1-\beta}{k_c+(d-1)k_j} .
\label{eq:coat}
\end{equation}
This relation defines the effective conductivity $k$ of the assemblage. Graded atoms, with a
core and a mantle whose composition depends on the radius, are built in the same way, layer by layer.

\section{The planar problem}\label{sec:2d}
\subsection{Bounds: derivation}\label{sec:bounds2d}

The planar problem was solved by Cherkaev [2009, 2012]. In this section, we present his results,
as they will be required to solve the three-dimensional problem. The bound has the form
\begin{equation}
2k_0\ \ge\ 2k_*=\max_{r,\tau}V_{2D}(r,\tau),\qquad V_{2D}(r,\tau)=\sum_{i=1,2}m_i\,g_i(r,\tau)+2r-2\tau ,
\label{eq:wd2d}
\end{equation}
\begin{equation}
g_i(r,\tau)=\min_{\det E\ge0}\bigl[k_i|E|^2+2\tau\det E-r\Tr E\bigr].
\label{eq:g2d0}
\end{equation}
We already mentioned (Section~\ref{sec:translation}) that
for $r>0$ the minimizer in \eqref{eq:g2d0} is a symmetric nonnegative
matrix with real nonnegative eigenvalues $\lambda_i$, and
\begin{equation}
g_i(r,\tau)=\min_{\lambda_1,\lambda_2\ge0}\bigl[k_i(\lambda_1+\lambda_2)^2+2(\tau-k_i)\lambda_1\lambda_2
-r(\lambda_1+\lambda_2)\bigr],
\label{eq:g2d}
\end{equation}
the minimum \eqref{eq:gilambda} over the cone $E\succeq0$, written in the eigenvalues of $E$. 
In terms of $T=\lambda_1+\lambda_2$ and $D=\lambda_1\lambda_2\in [0,T^2/4]$
the bracket in \eqref{eq:g2d} is
\begin{equation}
q_i=k_i T^2+2(\tau-k_i) D-r\,T.
\label{eq:rule2d}
\end{equation}
It is minimized over $\lambda_1,~\lambda_2$ with the fixed trace $T$;
the minimizer is
\newline $\bullet$ the isotropic field, $D=T^2/4$, if $\tau<k_i$;
\newline $\bullet$ a rank-one field, $D=0$, if $\tau>k_i$;
\newline $\bullet$ and any field of the minimizing trace if $\tau=k_i$; the Lagrangian is degenerate, $q_i=k_iT^2-rT$.
\newline
For $\tau<k_1$, both materials have isotropic minimizers and $V_{2D}(\tau)$ increases with $\tau$ up to $\tau=k_1$.
For $\tau>k_2$ both minimizers are rank one, and $V_{2D}$ decreases with $\tau$.
Therefore, $\tau\in[k_1,k_2]$. 

Inside the interval, $\tau\in (k_1,k_2)$,
the trace in $S_1$ is
$T=r/(2k_1)$, $D=0$ and $g_1=-r^2/(4k_1)$, the field being rank one.
 In $S_2$ the field is isotropic, $\lambda_1=\lambda_2=r/(2(k_2+\tau))$, and
$g_2=-r^2/(2(k_2+\tau))$.
The maximum over $r$ is at $r=2(k_*+\tau)$, and
\begin{equation}
k_*=\max_{k_1\le\tau\le k_2}\Psi(\tau),\qquad \Psi(\tau)=\Bigl[\frac{m_1}{2k_1}+\frac{m_2}{k_2+\tau}\Bigr]^{-1}-\tau .
\label{eq:F2d}
\end{equation}

The function $\Psi(\tau)$ is concave, and its
stationary point is
\begin{equation}
k_2+\tau=\frac{2k_1\bigl(\sqrt{m_2}-m_2\bigr)}{m_1} .
\label{eq:tau2d}
\end{equation}
It lies in $[k_1,k_2]$ when $m_{12}\le m_1\le m_{11}$, where
\begin{equation}
m_{11}=\frac{2k_1}{k_1+k_2}\bigl(\sqrt{m_2}-m_2\bigr),\qquad
m_{12}=\frac{k_1}{k_2}\bigl(\sqrt{m_2}-m_2\bigr),
\label{eq:m11}
\end{equation}
and otherwise the maximum is at the nearer end of the interval. 

At  $\tau=k_1$ any field of the same trace gives the same value of $g_1$. 
At  $\tau=k_2$ any field of the same trace gives the same value of $g_2$. 
This gives three regimes. In the order of decreasing $m_1$, we call them the
\emph{$S_1$-wheel regime}, which is the Hashin--Shtrikman regime, the \emph{bare-wheel regime} and
the \emph{$S_2$-wheel regime}, and denote the bound in them by $k_{01}$, $k_{02}$ and $k_{03}$. In
the $S_1$-wheel and $S_2$-wheel regimes the bound coincides with that of Nesi [1995]. The names
refer to the optimal structures and will become clear in Section~\ref{sec:structures2d}.

\paragraph{The $S_1$-wheel regime, $m_1\ge m_{11}$.} Here $\tau=k_1$, $r=2(k_{01}+k_1)$, and
\begin{equation}
k_{01}=\Bigl[\frac{m_1}{2k_1}+\frac{m_2}{k_1+k_2}\Bigr]^{-1}-k_1 .
\label{eq:k01}
\end{equation}
The Lagrangian of $S_1$ is degenerate, and only the trace of $E_1$ is fixed. The support sets are
\begin{equation}
\Lambda_1=\{\lambda\ge0:\ \lambda_1+\lambda_2=2\alpha_1\},\qquad \Lambda_2=\{\alpha_2(1,1)\},\qquad
\alpha_1=\frac{k_{01}+k_1}{2k_1},\quad \alpha_2=\frac{k_{01}+k_1}{k_1+k_2}.
\label{eq:sets1}
\end{equation}
The integral restriction on the trace reads $m_1\alpha_1+m_2\alpha_2=1$. The integral restriction
on the determinant fixes the average determinant of $S_1$, $m_1\langle\det E\rangle_1=1-m_2
\alpha_2^2$. On the segment $\Lambda_1$ the determinant runs from $0$ at its ends to $\alpha_1^2$
at its midpoint, so this restriction can be satisfied if and only if
\begin{equation}
m_2\alpha_2^2\le1\qquad(\text{equivalently } m_1\ge m_{11});
\label{eq:cond1}
\end{equation}
the other inequality, $m_1\alpha_1^2+m_2\alpha_2^2\ge1$, follows from the restriction on the trace.
At
$m_1=m_{11}$ the field $E_1$ takes only the two end values of its segment.

\paragraph{The bare-wheel regime, $m_{12}\le m_1\le m_{11}$.} Here $\tau$ is \eqref{eq:tau2d},
$r=4k_1(1-\sqrt{m_2})/m_1$, and
\begin{equation}
k_{02}=k_2+\frac{2k_1}{m_1}\bigl(1-\sqrt{m_2}\bigr)^2 .
\label{eq:k02}
\end{equation}
The field of $S_1$ is rank one and that of $S_2$ isotropic:
\begin{equation}
\Lambda_1=\{(c,0),(0,c)\},\quad c=\frac{2\bigl(1-\sqrt{m_2}\bigr)}{m_1},\qquad
\Lambda_2=\Bigl\{\frac{1}{\sqrt{m_2}}(1,1)\Bigr\}.
\label{eq:sets2}
\end{equation}
These values follow from the integral restrictions alone.
Since $\det E_1=0$, the restriction on the determinant gives $m_2\det E_2=1$,
which
defines the isotropic field of $S_2$, $E_2=I/\sqrt{m_2}$, independently of $m_1$ and of $\tau$; the
restriction on the trace, $m_1c+2\sqrt{m_2}=2$, then gives $c$. Stationarity of $V_{2D}$ in $\tau$ gives
the integral restriction on the determinant, and stationarity in $r$ gives the one on the trace.

\paragraph{The $S_2$-wheel regime, $m_1\le m_{12}$.} Here $\tau=k_2$, $r=2(k_{03}+k_2)$, and
\begin{equation}
k_{03}=\Bigl[\frac{m_1}{2k_1}+\frac{m_2}{2k_2}\Bigr]^{-1}-k_2 .
\label{eq:k03}
\end{equation}
Now the Lagrangian of $S_2$ is degenerate:
\begin{equation}
\Lambda_1=\{(c,0),(0,c)\},\quad c=\frac{2k_2}{k_1a_*},\qquad \Lambda_2=\Bigl\{\lambda\ge0:\ \lambda_1+\lambda_2=\frac{2}{a_*}\Bigr\},\qquad a_*=\frac{k_2}{k_1}m_1+m_2 .
\label{eq:sets3}
\end{equation}
The integral restriction on the determinant falls on $S_2$ alone, $m_2\langle\det E\rangle_2=1$, and
on its segment the determinant is at most $1/a_*^2$; so the $S_2$-wheel regime is
\begin{equation}
a_*^2\le m_2\qquad(\text{equivalently } m_1\le m_{12}).
\label{eq:cond3}
\end{equation}
At $m_1=m_{12}$ the support of $S_2$ collapses to the midpoint of its segment, the isotropic point of the bare-wheel regime.

\subsection{Bounds: illustration}\label{sec:bounds2dill}

\begin{figure}[htbp]
	\begin{center}
		\IfFileExists{nov18_bound2d.png}{\includegraphics[width=0.90\linewidth]{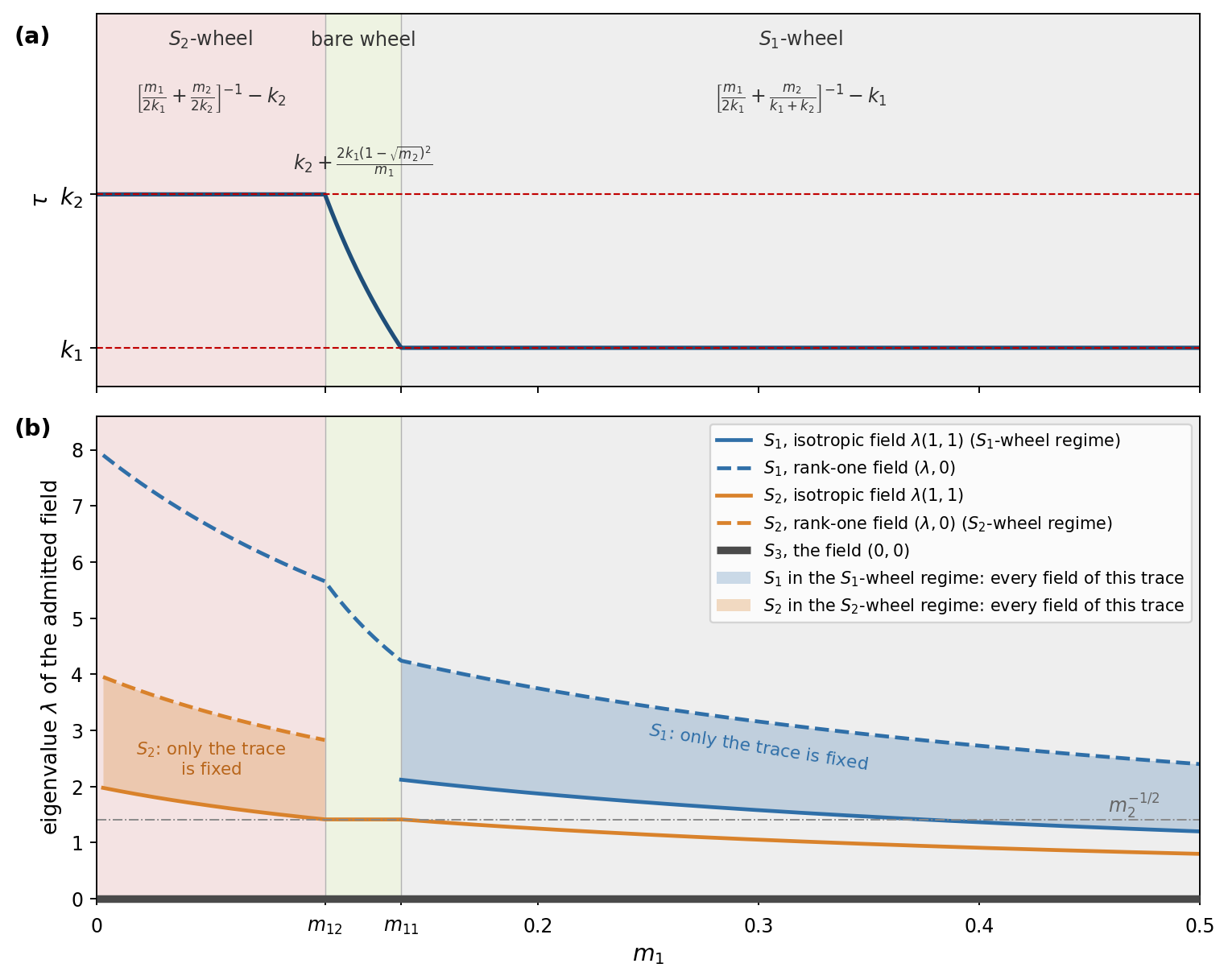}}{%
			\framebox[0.8\linewidth][c]{\rule{0pt}{4cm}\textit{[figure file not found]}}}
		\caption{The planar bound, $k_1=1$, $k_2=2$, $m_2=\tfrac12$: (a) the multiplier $\tau$ and the
			bound in each regime; (b) the fields admitted in $S_1$, $S_2$ and $S_3$.}
		\label{fig:bound2d}
	\end{center}
\end{figure}

The sequence of the regimes is illustrated in Figure~\ref{fig:bound2d}.  As $m_1$ decreases, $\tau$ moves from $k_1$ to
$k_2$ (Figure~\ref{fig:bound2d}(a)). $S_1$ passes from a field with only
the trace fixed to a rank-one field at $m_{11}$, and $S_2$ passes from an isotropic field to a field
with only the trace fixed at $m_{12}$. The bound and its derivative in
$m_1$ are continuous at both thresholds, and so are the field values.
At $m_{11}$ the ends of $\Lambda_1$ are the rank-one points of the bare-wheel regime, and at $m_{12}$ the midpoint
of $\Lambda_2$ is the isotropic point of the bare-wheel regime.

Figure~\ref{fig:bound2d}(b) shows the admitted fields by their eigenvalue $\lambda$: dashed lines
for the rank-one fields $(\lambda,0)$, solid lines for the isotropic fields $\lambda(1,1)$, and the
dark line for the field $(0,0)$ of $S_3$. In the bare-wheel and $S_2$-wheel regimes the field of $S_1$ is rank one, with
$\lambda=c$; in the $S_1$-wheel regime every field of $S_1$ of the trace $2\alpha_1$ is admitted, from the
rank-one field $(2\alpha_1,0)$ to the isotropic field $\alpha_1(1,1)$, and the shading spans the
larger eigenvalues of these fields. Likewise, the field of $S_2$ is isotropic in the $S_1$-wheel and bare-wheel regimes,
with $\lambda=\alpha_2$, and $\lambda=m_2^{-1/2}$ (dash-dot line) independently of $m_1$ in the bare-wheel regime; in the $S_2$-wheel regime every field of $S_2$ of the trace $2/a_*$ is admitted, from $(1,1)/a_*$ to $(2/a_*,0)$.

\begin{figure}[htbp]
\begin{center}
\IfFileExists{nov18b_supportsets2d.png}{\includegraphics[width=0.90\linewidth]{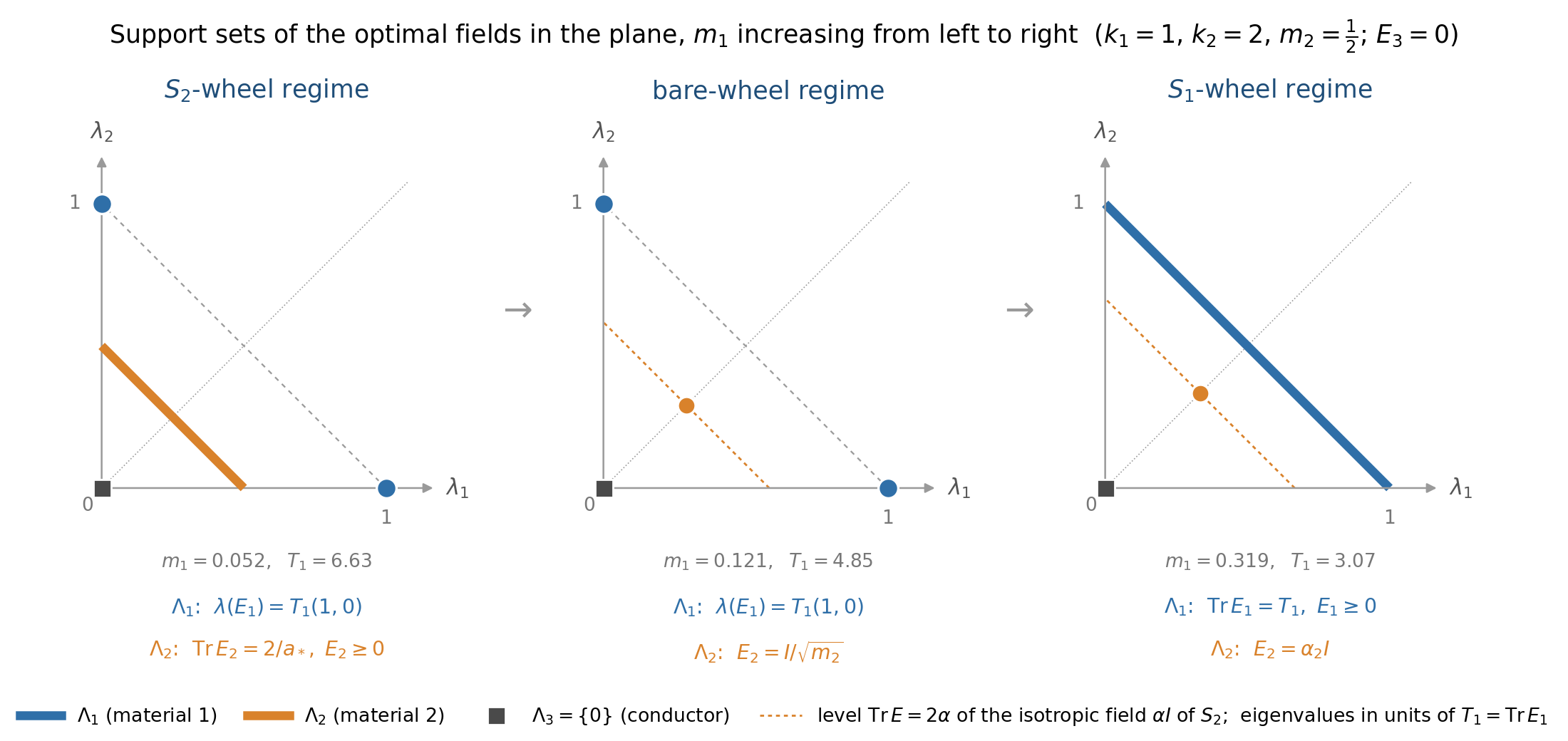}}{%
\framebox[0.8\linewidth][c]{\rule{0pt}{4cm}\textit{[figure file not found]}}}
\caption{The support sets of the planar bound in the three regimes, $k_1=1$, $k_2=2$,
$m_2=\tfrac12$.}
\label{fig:sets2d}
\end{center}
\end{figure}

The support sets are illustrated in Figure~\ref{fig:sets2d}, for one value of $m_1$ in the middle
of each regime. In every regime $T_1=\Tr E_1=r/(2k_1)$ on $\Lambda_1$; the figure shows the eigenvalues in units of $T_1$, and the dashed orange line marks the trace of the isotropic field of $S_2$. Thus $|E_1|\le T_1$, with equality at the rank-one fields. In the $S_1$-wheel regime, $|E_1|\ge\sqrt2\,\alpha_1$ on the segment $\Lambda_1$ and
$|E_2|=\sqrt2\,\alpha_2$, so $|E_2|/|E_1|\le\alpha_2/\alpha_1=2k_1/(k_1+k_2)$, with equality at the
midpoint. In the bare-wheel regime, $|E_2|/|E_1|=\sqrt2/(\sqrt{m_2}\,c)=m_1/\bigl(\sqrt2(\sqrt{m_2}-m_2)\bigr)$,
which is largest at $m_{11}$, where it equals $\sqrt2\,k_1/(k_1+k_2)$. In the $S_2$-wheel regime, $|E_2|\le2/a_*$
and $|E_1|=c$, so $|E_2|/|E_1|\le k_1/k_2$. Thus $|E_2|\le\frac{2k_1}{k_1+k_2}|E_1|$ for every pair
of points of $S_1$ and $S_2$ in every regime, and $E_3=0$.

\subsection{Optimal structures}\label{sec:structures2d}

\begin{figure}[htbp]
	\begin{center}
		\IfFileExists{nov18_atoms2d.png}{\includegraphics[width=\linewidth]{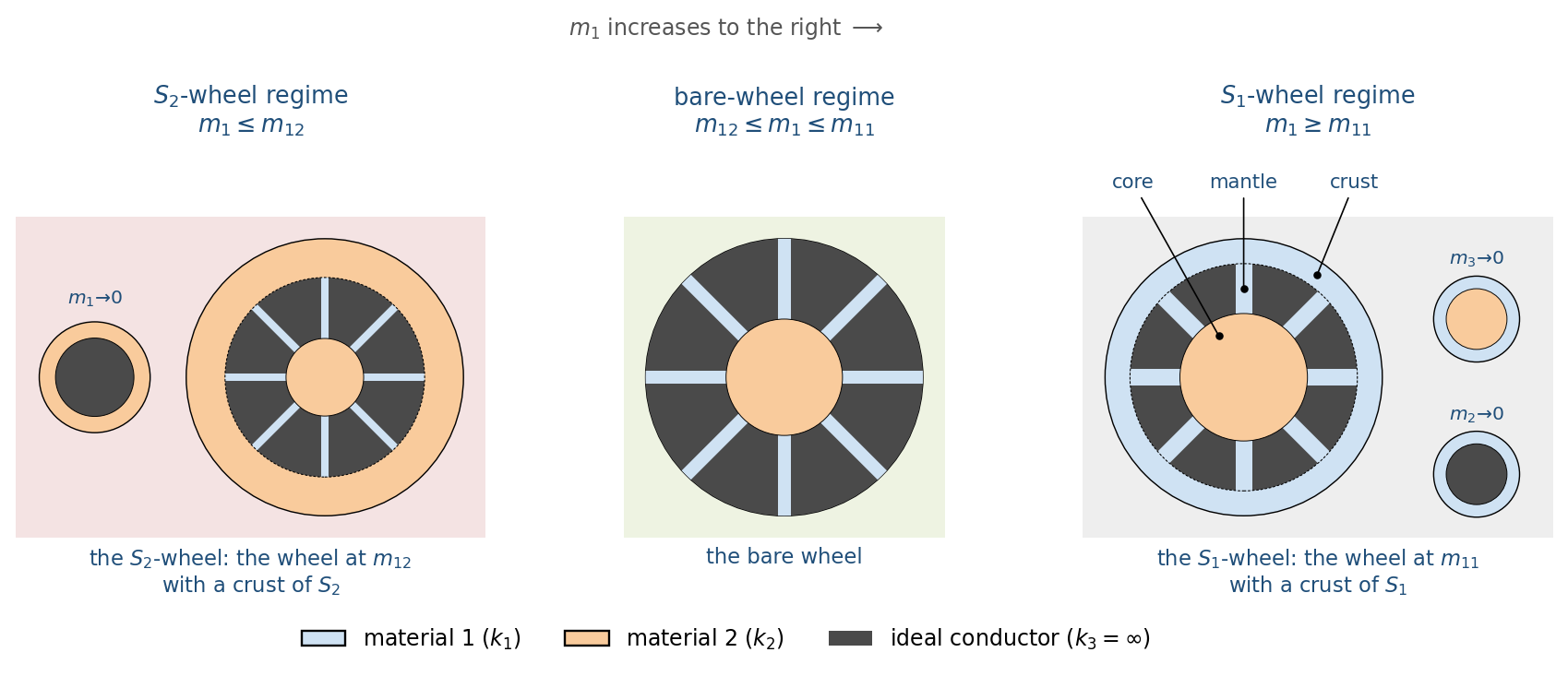}}{%
			\framebox[0.8\linewidth][c]{\rule{0pt}{4cm}\textit{[figure file not found]}}}
		\caption{The optimal atoms in the plane, $m_1$ increasing to the right.}
		\label{fig:atoms2d}
	\end{center}
\end{figure}

Any structure that realizes the support fields is optimal. In this section, we present one type of such structures, well
suited to describe isotropic materials, the wheels [Cherkaev 2012]. They are illustrated in
Figure~\ref{fig:atoms2d}, with their layers named: the bare wheel (center), a core of $S_2$ in a mantle of radial spikes of
$S_1$ and wedges of $S_3$, and the wheels at the thresholds $m_{12}$ and $m_{11}$ with a crust of
$S_2$ (left) and of $S_1$ (right), the $S_2$-wheel and the $S_1$-wheel, given by the enveloping
principle (Section~\ref{sec:envelope2d}). These three structures explain the names of the regimes. The small discs are the two-material limits: the disc of
$S_3$ coated by $S_2$ as $m_1\to0$, and the discs of $S_2$ and of $S_3$ coated by $S_1$ as $m_3\to0$
and as $m_2\to0$. The spikes of $S_1$ have a constant width, and the wedges of $S_3$ between them
meet the core along arcs \eqref{eq:theta2d}.

\subsubsection{The bare-wheel regime: the wheel}\label{sec:wheel2d}

We start with a structure that attains the bound in the bare-wheel regime, the wheel assemblage; without a crust, we call it the bare wheel.
In this regime, the bound prescribes a rank-one field in $S_1$ with the fixed nonzero eigenvalue $c$, the
isotropic field $I/\sqrt{m_2}$ in $S_2$, and $E=0$ in $S_3$ \eqref{eq:sets2}. A rank-one
field is the field of a laminate with $S_3$. Across a flat interface, the tangential components of each column
of $E$ are continuous, and on the $S_3$ side they vanish. In a layer of $S_1$ between two faces
of the conductor with the normal $n$, every column of $E$ is therefore parallel to $n$,
$E=n\otimes v$, and $E\succeq0$ forces $E=\lambda\,n\otimes n$.

A laminate of $S_1$ and $S_3$ with the normal $n$ and the fraction $\theta$ of $S_1$ has the
conductivity $k_1/\theta$ across the layers and an infinite conductivity along them. Physically, the
current runs freely along the layers of the conductor, but across them it must pass the layers of
$S_1$ in series.

To obtain an isotropic structure, let us coil the laminate around a central disc of $S_2$, the
core. The layers become radial spikes. The $S_1$ spikes have a constant width, while the conductor $S_3$ wedges
between them widen to fill the mantle. The mantle, the annulus around the core, conducts infinitely
along the radius and with $k_1/\theta(\rho)$ around the circle. The radial conductor carries the potential of the outer circle
to the disc unchanged. Indeed, a unit field outside the unit circle imposes there the potential
$\cos\varphi$. The same potential acts on the circle of radius $\sqrt{m_2}$ that bounds the disc, and the field
in the disc is $(1/\sqrt{m_2}) I$. It is isotropic, as the bound requires.
The current's path is split: one part flows through the core, and another flows around it through the spikes of $S_1$.

At the radius $\rho$ the tangential field, averaged over the annulus, is $1/\rho$.  The
field in a spike is therefore $(\rho\,\theta(\rho))^{-1}e_\varphi\otimes e_\varphi$. The bound fixes
its eigenvalue at $c$, and we obtain
\begin{equation}
\theta(\rho)=\frac{1}{c\,\rho},\qquad \sqrt{m_2}\le\rho\le1 .
\label{eq:theta2d}
\end{equation}
The volume restriction of $S_1$,
$2\int_{\sqrt{m_2}}^1\theta\rho\,d\rho=2(1-\sqrt{m_2})/c=m_1$, holds identically. Besides,
$\theta\le1$: its largest value, at $\rho=\sqrt{m_2}$, is $m_1/(2(\sqrt{m_2}-m_2))$, which runs from
$k_1/(2k_2)$ at $m_{12}$ to $k_1/(k_1+k_2)$ at $m_{11}$.
Once the laminate is coiled around a core, the bound fixes its grading.
The wheel attains the bound since all fields of the wheel lie in the support sets.

\subsubsection{The \texorpdfstring{$S_1$}{S1}-wheel and \texorpdfstring{$S_2$}{S2}-wheel regimes: the enveloping principle}\label{sec:envelope2d}

Here we show that the outer regimes need not be treated separately. We will use the enveloping principle, which is applicable due to a special algebraic form of the bounds in these two regimes. Note that $S_1$ enters the
bound of the $S_1$-wheel regime, and $S_2$ the bound of the $S_2$-wheel regime, exactly as a coat enters the conductivity of a
coated disc.

Consider an assemblage of discs of any structure that attains the bound in the $S_1$-wheel regime, each coated by $S_1$.
Let us show that the resulting structure also attains the bound at the new
fractions. In the $S_2$-wheel regime the same is true of a coat of $S_2$.

As mentioned in Section~\ref{sec:structures}, a disc (sphere) of conductivity $k_c$ coated by an annulus (shell) of $S_j$ that takes the fraction $\beta$ of its volume has the effective conductivity $k$ given by \eqref{eq:coat}. If that inclusion is coated with an additional amount $\hat{\beta}$ of $S_j$, its effective conductivity satisfies \eqref{eq:coat} with the adjusted fraction $\beta\to\hat{\beta}+(1-\hat{\beta})\beta$: two coats act as one.

A coat of $S_1$ that takes the fraction $\beta$ maps the fractions
$\hat m=(\hat m_1,\hat m_2)$ of the core to
\[
m_1=\beta+(1-\beta)\hat m_1,\qquad m_2=(1-\beta)\hat m_2 .
\]
As $\beta$ grows from $0$ to $1$, the point $m$ runs along the straight segment from $\hat m$ to the
vertex $m_1=1$ of pure $S_1$; along this ray the composition of the core, the ratio $m_2/(m_2+m_3)$, does not change. Likewise, a coat of $S_2$ moves $m$ along the segment from
$\hat m$ to the vertex $m_2=1$ of pure $S_2$, and the ratio $m_1/(m_1+m_3)$ does not change.

The bounds \eqref{eq:k01} and \eqref{eq:k03} of the $S_1$-wheel and $S_2$-wheel regimes read
\[
\frac{1}{k_{01}+k_1}=\frac{m_1}{2k_1}+\frac{m_2}{k_1+k_2},\qquad
\frac{1}{k_{03}+k_2}=\frac{m_1}{2k_1}+\frac{m_2}{2k_2}.
\]
Both are affine in the fractions, and the weight of $S_1$ in the first, and of $S_2$ in the second,
is $1/(2k_j)$, the weight of a coat of $S_j$ in \eqref{eq:coat} with $d=2$. Substituting the fractions of the ray into the first bound turns it into \eqref{eq:coat} with $k_c=k_{01}(\hat m)$, and likewise for the second bound with $k_c=k_{03}(\hat m)$. Therefore, if a structure
attains the bound of the $S_1$-wheel regime at $\hat m$, the assemblage of its discs coated by $S_1$ attains it at
every point of the ray toward pure $S_1$ that lies in the $S_1$-wheel regime; the same holds for a coat of $S_2$ in the $S_2$-wheel regime. Since $2k_{01}=\max_rV_{2D}(r,k_1)$ and $2k_{03}=\max_rV_{2D}(r,k_2)$ are values of the bound
functional at every fraction, the coated structures are optimal.

\begin{figure}[htbp]
 	\begin{center}
 		\IfFileExists{nov18b_regions2d.png}{\includegraphics[width=0.95\linewidth]{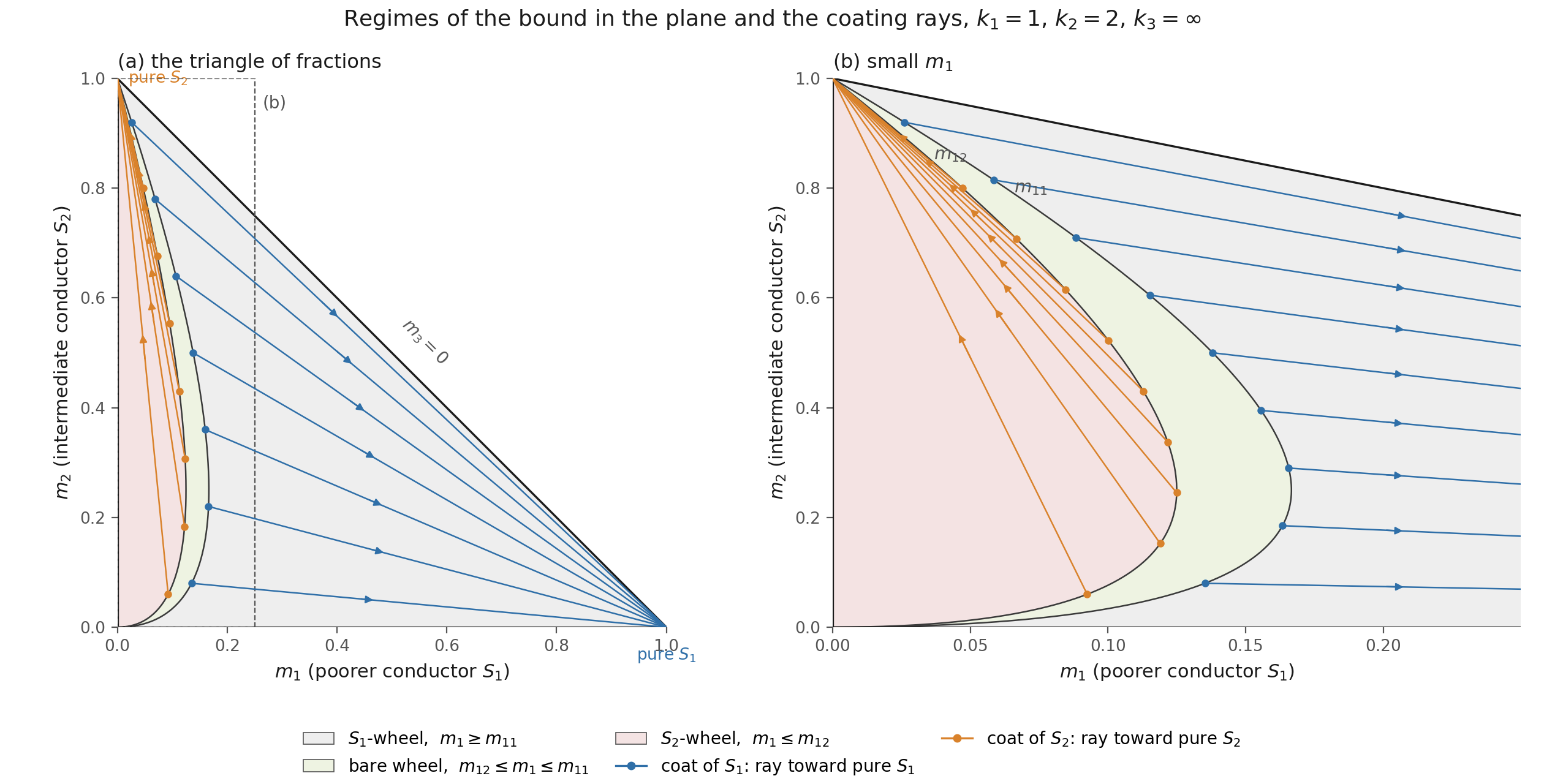}}{%
 			\framebox[0.8\linewidth][c]{\rule{0pt}{4cm}\textit{[figure file not found]}}}
 		\caption{The three planar regimes and the coating rays in the plane of the fractions
 			$(m_1,m_2)$, $k_1=1$, $k_2=2$: (a) the whole triangle, (b) the strip of small $m_1$.}
 		\label{fig:regimes2d}
 	\end{center}
 \end{figure}

The rays that start on the threshold curve
$m_1=m_{11}(m_2)$ and end at the vertex of pure $S_1$ stay in the $S_1$-wheel regime and cover it, each point of
the regime lying on exactly one ray; the rays from the curve $m_1=m_{12}(m_2)$ to the vertex of pure
$S_2$ cover the $S_2$-wheel regime in the same way (Figure~\ref{fig:regimes2d}). Therefore the structures at the
two boundaries of the bare-wheel regime, coated by $S_1$ or by $S_2$, are optimal in the whole of the $S_1$-wheel and $S_2$-wheel regimes, and new structures need to be constructed only in the bare-wheel regime.

For the coat of $S_1$ in the $S_1$-wheel regime this argument is due to Cherkaev [2012, \S4.2, eq.~(33)], where
the details can be found; Cherkaev [2009, App.~A.2] gives its laminate form, and Albin and Cherkaev
[2006, \S3.2] coat their structures with layers of $S_1$ in space. The coat of $S_2$ in the $S_2$-wheel regime follows in the same way, with $k_2$ in place of $k_1$.

In Figure~\ref{fig:regimes2d}, a structure that attains the bound at a point of the curve
$m_1=m_{11}$ (blue dot), coated by $S_1$, attains it along the blue ray from this point to the
vertex of pure $S_1$; a structure at a point of the curve $m_1=m_{12}$ (orange dot), coated by
$S_2$, attains it along the orange ray to the vertex of pure $S_2$. The arrows show the direction
in which the fractions move as the coat grows.

Thus the $S_1$-wheel and the $S_2$-wheel, the wheel coated by $S_1$ or by $S_2$, are optimal in the outer regimes.

\section{The three-dimensional problem}\label{sec:3d}

Here we apply the described method to three-dimensional composites.
\subsection{The bound, region by region}\label{sec:bound3d}

 With $d=3$ the bound \eqref{eq:wd}--\eqref{eq:gilambda} reads as follows:
\begin{equation}
3k=\sum_{i=1,2}m_i\Bigl\langle k_i|E|^2-r\Tr E+\tau\Phi(E)+\mu\det E\Bigr\rangle_i+3r-6\tau-\mu
\ \ge\ V(r,\tau,\mu),
\label{eq:wd3d}
\end{equation}
\begin{equation}
V(r,\tau,\mu)=m_1g_1(r,\tau,\mu)+m_2g_2(r,\tau,\mu)+3r-6\tau-\mu,
\label{eq:V}
\end{equation}
\begin{equation}
g_i(r,\tau,\mu)=\min_{\lambda\ge0}\Bigl[k_i\sum_j\lambda_j^2-r\sum_j\lambda_j+\tau\Phi(\lambda)
+\mu\prod_j\lambda_j\Bigr],
\label{eq:g}
\end{equation}
where $\lambda$ are the eigenvalues of $E$. 
As mentioned above (Section~\ref{sec:translation}), a minimizer $E$ is
symmetric and its eigenvalues $\lambda_i$ are real and (by the restriction \eqref{eq:psd}) nonnegative,
 $\lambda_i\geq 0$.

The bound is
\begin{equation}
3k_*=\max_{r,\ \tau,\ \mu\ge0}V(r,\tau,\mu) . \label{eq:bound3d}
\end{equation}

\begin{figure}[htbp]
 	\centering
 	\IfFileExists{nov18b_multipliers3d.png}{\includegraphics[width=0.9\linewidth]{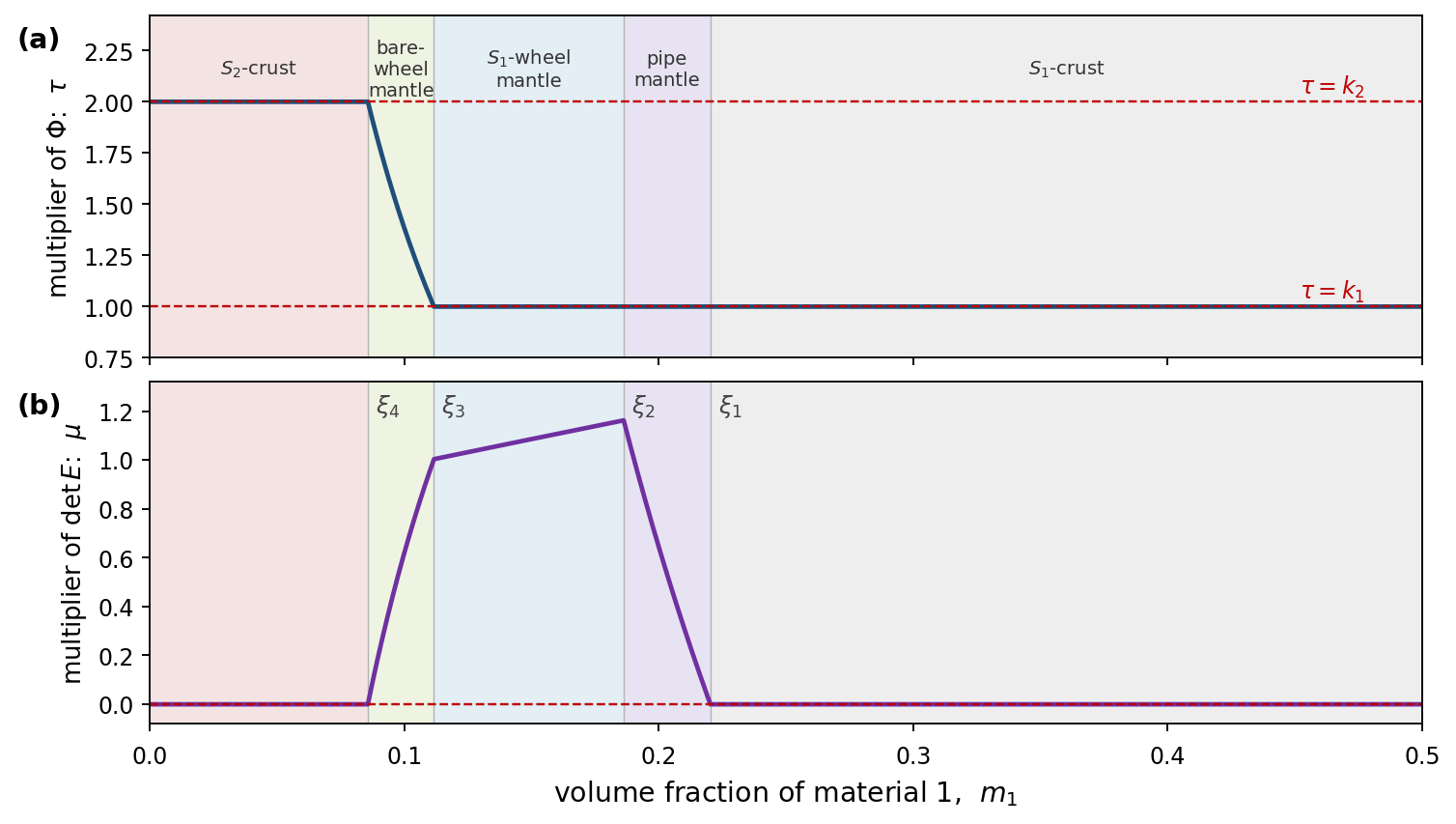}}{%
 		\framebox[0.8\linewidth][c]{\rule{0pt}{4cm}\textit{[figure file not found]}}}
 	\caption{The multipliers (a) $\tau$ and (b) $\mu$ in the five regions of the three-dimensional
 		bound, $k_1=1$, $k_2=2$, $m_2=\tfrac12$; thin grey vertical lines show their boundaries $\xi_1,\ldots,\xi_4$.}
 	\label{fig:multipliers3d}
 \end{figure}

The regions are listed in Table~\ref{tab:regions}, and the multipliers are shown in Figure~\ref{fig:multipliers3d}.
In the order of decreasing $m_1$, we call the five regions the \emph{$S_1$-crust region}, which
is the Hashin--Shtrikman region, the \emph{pipe mantle}, \emph{$S_1$-wheel mantle} and
\emph{bare-wheel mantle} regions, and the \emph{$S_2$-crust region}. The names refer to the optimal
structures and will become clear in Section~\ref{sec:structures3d}. In the following, we show the
derivation of the bound in each of the five regions.

 \begin{table}[htbp]
 	\centering
 	\begin{tabular}{@{}lllll@{}}
 		\hline
 		range of $m_1$ & region & bound & $\tau$ & $\mu$\\
 		\hline
 		$[\xi_1,1-m_2]$ & $S_1$-crust & $B^{(1)}$ \eqref{eq:B1} & $k_1$ & $0$\\
 		$[\xi_2,\xi_1]$ & pipe mantle & $B_2^{*}$ \eqref{eq:B2} & $k_1$ & $>0$\\
 		$[\xi_3,\xi_2]$ & $S_1$-wheel mantle & cubic \eqref{eq:cubic}, \eqref{eq:wsbound} & $k_1$ & $>0$\\
 		$(\xi_4,\xi_3)$ & bare-wheel mantle & pair \eqref{eq:pair}, \eqref{eq:wsbound} & in $(k_1,k_2)$ & $>0$\\
 		$(0,\xi_4]$ & $S_2$-crust & $B^{(3)}$ \eqref{eq:B3} & $k_2$ & $0$\\
 		\hline
 	\end{tabular}
 	\caption{The regions of the three-dimensional bound.}
 	\label{tab:regions}
 \end{table}

 \subsection{The bound, derivation}\label{sec:bound3dder}

In terms of the invariants $T=\sum\lambda_j$, $\Phi$ and
$D=\prod\lambda_j$ the bracket in \eqref{eq:g} becomes
\begin{equation}
q_i=k_iT^2+(\tau-k_i)\Phi+\mu D-rT .
\label{eq:q3d}
\end{equation}
As in the plane problem, $\tau<k_i$ leans toward the isotropic field, $\tau>k_i$ sends the field to the boundary of the octant, and at $\tau=k_i$ the second invariant drops out. The Lagrange multiplier $\mu>0$
penalizes the determinant. 

Throughout what follows we write
\[
A=k_2+2\tau,\qquad C=k_2+\tau,\qquad \kappa=\frac{1-\sqrt{1-A/(2C)}}{A},\qquad a=m_2^{1/3}.
\]
Let us compute the densities of two fields of $S_2$: the isotropic field $(b,b,b)$ has the density
$3Ab^2+\mu b^3-3rb$, and the rank-two field $(\gamma,\gamma,0)$ has $2C\gamma^2-2r\gamma$, minimal at
$\gamma=r/(2C)$. Four thresholds of $m_1$, $\xi_1>\xi_2>\xi_3>\xi_4$, numbered in the order in which they are met as $m_1$ decreases, separate five regions (Table~\ref{tab:regions}). In all regions material $S_1$ has $\Tr E_1=r/(2k_1)$ and
$g_1=-r^2/(4k_1)$, and $E_3=0$.

\begin{remark}[The minimizers in \eqref{eq:g}]
For $S_1$ with $\tau\ge k_1$ and $\mu\ge0$, the terms $(\tau-k_1)\Phi$ and $\mu D$ in
\eqref{eq:q3d} are nonnegative on the octant $\lambda\ge0$. The minimizers therefore have the
trace $r/(2k_1)$ and make these terms vanish: any field of this trace if $\tau=k_1$ and $\mu=0$,
the fields with $D=0$ if $\tau=k_1$ and $\mu>0$, and the rank-one field if $\tau>k_1$. For $S_2$
with $\tau<k_2$, the bracket $q_2=k_2\sum_j\lambda_j^2+\tau\Phi+\mu D-rT$ grows at infinity, so its
minimum over the octant is attained at a critical point, in the interior or on the boundary. In the
interior, subtracting the equations $\partial q_2/\partial\lambda_i=0$ and
$\partial q_2/\partial\lambda_j=0$ gives $(\lambda_i-\lambda_j)\bigl(2(k_2-\tau)-\mu\lambda_l\bigr)=0$,
$\{i,j,l\}=\{1,2,3\}$. A critical point is therefore isotropic or, if $\mu>0$, has two eigenvalues
equal to $x=2(k_2-\tau)/\mu$ and the third $y\ne x$. At such a point the Hessian of $q_2$ has the
eigenvalue $2(k_2-\tau)-\mu y$ on $(1,-1,0)$, and its determinant on the orthogonal plane is
$2k_2\bigl(\mu y-2(k_2-\tau)\bigr)$; one of the two is negative, and the point is a saddle. On the
face $\lambda_3=0$, $q_2$ is a convex quadratic with the minimizer $(\gamma,\gamma,0)$,
$\gamma=r/(2C)$; on an edge its minimum $-r^2/(4k_2)$ exceeds $-r^2/(2C)$. Hence $g_2$ is attained
at the isotropic field or at the rank-two field $(\gamma,\gamma,0)$, the two fields compared above.
For $\tau=k_2$ and $\mu=0$ the Lagrangian of $S_2$ is degenerate, and every field of the trace
$r/(2k_2)$ is a minimizer.
\end{remark}

\paragraph{The $S_1$-crust region, $m_1\in[\xi_1,1-m_2]$.} Here
$$\tau=k_1, \quad \mu=0, \quad
r=2(B^{(1)}+2k_1),$$
and the bound is the classical translation (Hashin--Shtrikman) bound
\begin{equation}
B^{(1)}=\Bigl[\frac{m_1}{3k_1}+\frac{m_2}{k_2+2k_1}\Bigr]^{-1}-2k_1 .
\label{eq:B1}
\end{equation}
Only the trace of $E_1$ is fixed, and $E_2$ is isotropic,
\begin{eqnarray}
\Lambda_1=\Bigl\{\lambda\ge0:\ \sum_j\lambda_j=3\alpha_1\Bigr\},\qquad
\Lambda_2=\{b(1,1,1)\}, \label{eq:fields-hs} \\ b=\frac{B^{(1)}+2k_1}{k_2+2k_1},\qquad \frac{\alpha_1}{b}=\frac{k_2+2k_1}{3k_1}.
\nonumber
\end{eqnarray}

The determinant $D$ in $S_1$ varies; it reaches zero on the boundary of $\Lambda_1$. In $S_2$, $D=b^3$. The integral restriction $\langle D \rangle =1$ on the
determinant can be satisfied if and only if
\begin{equation}
m_2b^3\le1,\qquad\text{that is,}\qquad m_1\ge\xi_1=\frac{3k_1\bigl(m_2^{1/3}-m_2\bigr)}{k_2+2k_1}.
\label{eq:xi1}
\end{equation}
The integral restriction on $\Phi$ is met as well: the structures of
Section~\ref{sec:envelopes3d} attain $B^{(1)}$ in the whole region, and their fields satisfy all
three restrictions.

\begin{remark}
The necessity of the threshold
$\xi_1$ is the null Lagrangian argument of Liu [2011], which our derivation repeats in the
language of the integral restriction on $\det E$. For laminates it was shown earlier by Briane and
Nesi [2004, Th.~1.2], whose condition becomes $m_1<\xi_1$ for an ideal conductor. The sufficiency
of the threshold, the attainment of $B^{(1)}$ down to $\xi_1$, is due to the hairy spheres of
Albin and Cherkaev [2006], whose threshold becomes $\xi_1$ as $k_3\to\infty$. Below $\xi_1$ the
Hashin--Shtrikman fields can no longer satisfy the integral restriction on the determinant, which
then requires $\mu>0$, and the bound changes its form.
\end{remark}

\paragraph{The pipe mantle region, $m_1\in[\xi_2,\xi_1]$.} Here
$$\tau=k_1, \quad r=6k_1(1-a^2)/m_1, \quad
\mu=ra^2-2(k_2+2k_1)a>0.$$
In this region, $\det E_1=0$ in $S_1$ (the edges of the triangle $\sum_j\lambda_j=3\alpha_1$), and
$\Tr E_1=3(1-a^2)/m_1$. The integral restriction on the determinant then rests on $S_2$ alone,
$m_2\det E_2=1$, and it fixes the field of $S_2$ at $E_2=I/a$. The fields of the two materials thus
satisfy
\begin{eqnarray}
\Lambda_1=\Bigl\{\lambda\ge0:\ \sum_j\lambda_j=3\alpha_1, \quad    \prod_j\lambda_j =0\Bigr\},\qquad
\Lambda_2=\{b(1,1,1)\}, \label{eq:fields-sb}\\ b=m_2^{-1/3}=\frac1a,\qquad \qquad 6k_1\alpha_1=2(2k_1+k_2)b+\mu\,b^{2}.
 \nonumber
\end{eqnarray}
%
The bound is
\begin{equation}
B_2^{*}=k_2a+k_1\Bigl[\frac{3(1-a^2)^2}{m_1}-2(1-a)\Bigr],\qquad a=m_2^{1/3}.
\label{eq:B2}
\end{equation}

\begin{remark}[Derivation]
Let us put $\tau=k_1$; the hedgehog sphere of Section~\ref{sec:spiky} attains the resulting value,
which confirms this choice (Section~\ref{sec:structures}). Since $\det E_1=0$, the part of \eqref{eq:g} for $S_1$ is a parabola in the trace, with the minimum
$g_1=-r^2/(4k_1)$. The field in $S_2$ is $\lambda=(b,b,b)$ with $b=1/a$, so that
$g_2=3Ab^2-3rb+\mu b^3$ and $m_2b^3=1$. In \eqref{eq:V} the term $m_2\mu b^3$ cancels $-\mu$, but
the requirement $\mu\ge0$ remains. Using $m_2b^2=a$ and $m_2b=a^2$, we are left with a quadratic
function of $r$,
\[
V=-\frac{m_1r^2}{4k_1}+3r(1-a^2)+3Aa-6k_1 ,
\]
whose maximum, at $r=6k_1(1-a^2)/m_1$, is
\[
V=3\bigl[\,3k_1(1-a^2)^2/m_1+(k_2+2k_1)a-2k_1\,\bigr]=3B_2^{*}.
\]
The multiplier $\mu=ra^2-2(k_2+2k_1)a$ follows from the stationarity of $g_2$ at $b=1/a$, and the
condition $\mu\ge0$ reads $r\ge2(k_2+2k_1)/a$, that is, $m_1\le\xi_1$.
\end{remark}

At $m_1=\xi_1$ the bound $B_2^*$ meets $B^{(1)}$ tangentially. The lower end $\xi_2$ is where the
isotropic point ceases to be the global minimizer of $g_2$. The rank-two field $(\gamma,\gamma,0)$
with $\gamma=r/(2C)$ gives $q_2=-r^2/(2C)$, which equals the isotropic value
$q_2= A/a^2-2r/a$ at $r=(2C+\sqrt{2Ck_2})/a$. In closed form,
\begin{equation}
\frac{\xi_2}{\xi_1}=\frac{2(k_2+2k_1)}{2(k_1+k_2)+\sqrt{2k_2(k_1+k_2)}},
\label{eq:xi2}
\end{equation}
free of $m_2$. The multiplier $\mu$ is largest there, $\mu_{\max}=m_2^{1/3}\bigl(\sqrt{2k_2(k_1+k_2)}
-2k_1\bigr)$, and its graph has a corner.

\paragraph{The $S_1$-wheel mantle region, $m_1\in[\xi_3,\xi_2]$.} Here $\tau=k_1$;  $r$ is the
positive root of the cubic
\begin{equation}
P\,r^3-3r^2+Q=0,\qquad P=\frac{m_1}{2k_1}+\frac{m_2}{C},\qquad Q=\frac{3}{\kappa^2}-\frac{1}{C\kappa^3},
\label{eq:cubic}
\end{equation}
and this root is unique because $Q<0$. With $\kappa=1/\bigl(2C(1+\sqrt{1-A/(2C)})\bigr)$ the
condition $Q<0$ reads $3C\kappa<1$, which holds for $\tau<k_2$. Once $r$ is found, the bound is the
explicit function
\begin{equation}
3k_*=-\frac{m_1r^2}{4k_1}-\frac{m_2r^2}{2C}+3r-6\tau+\frac{2A\kappa-1}{\kappa^2r},
\label{eq:wsbound}
\end{equation}
and
\begin{equation}
\mu=\frac{1-2A\kappa}{\kappa^2r}.
\label{eq:mu-window}
\end{equation}
Thus one root of a cubic and one substitution give the bound.

Optimal fields are:
 $$\Lambda_1:\quad  \det E_1=0,\qquad \Tr E_1=r/(2k_1).$$
 $S_2$ carries two field values,
\begin{equation}
\Lambda_2:\quad p\,\delta_{(b,b,b)}+(1-p)\,\delta_{(\gamma,\gamma,0)},\qquad
b=\kappa r,\quad \gamma=\frac{r}{2C},\quad p=\frac{1}{m_2b^3};
\label{eq:twoatom}
\end{equation}
the two have the same translated energy, and their proportion $p$ is fixed by the integral restriction
on the determinant.

\paragraph{The bare-wheel mantle region, $m_1\in(\xi_4,\xi_3)$.} Here $\tau > k_1$, and the
multipliers $r$ and $\tau$ are determined together by the stationarity of $V$ as the solutions to a
$2\times2$ algebraic system
\begin{equation}
\frac{m_1r}{2k_1}+\frac{3}{\kappa^2r^2}+\frac{m_2r}{C}-\frac{1}{C\kappa^3r^2}=3,\qquad
\frac{6}{\kappa r}+\frac{m_2r^2}{2C^2}-\frac{1}{2C^2\kappa^3r}=6,
\label{eq:pair}
\end{equation} 
The bound is again \eqref{eq:wsbound}, and $\mu$ is \eqref{eq:mu-window}. 

The support sets are
$$\Lambda_1=\bigl\{(r/(2k_1),0,0)\bigr\},\qquad \Lambda_2\ \ \text{as in \eqref{eq:twoatom}},$$
with $A$, $C$ and $\kappa$ at the current $\tau$. At the upper end $\xi_3$ the multiplier $\tau$
reaches $k_1$. The second equation of \eqref{eq:pair} does not contain $m_1$; multiplied by
$2C^2r$, at $\tau=k_1$ it becomes the cubic in \eqref{eq:xi3}, and $r$ is its largest positive
root, the only one with $p<1$ in \eqref{eq:twoatom} (Remark~\ref{rem:order}). The first equation
of \eqref{eq:pair} then gives $m_1=\xi_3$:
\begin{equation}
m_2r^{3}-12C^{2}r+2C^{2}\Bigl(\frac{6}{\kappa}-\frac{1}{2C^{2}\kappa^{3}}\Bigr)=0
\ \ \text{at }\tau=k_1,\qquad
\xi_3=\frac{2k_1}{r}\Bigl[3-\frac{3}{\kappa^{2}r^{2}}-\frac{m_2r}{C}+\frac{1}{C\kappa^{3}r^{2}}\Bigr].
\label{eq:xi3}
\end{equation}

\paragraph{The $S_2$-crust region, $m_1\in(0,\xi_4]$.} Here $\tau=k_2$, $\mu=0$ and
$r=2(B^{(3)}+2k_2)$. The bound again has the Hashin--Shtrikman form, with $\tau=k_2$,
\begin{equation}
B^{(3)}=\Bigl[\frac{m_1}{3k_1}+\frac{m_2}{3k_2}\Bigr]^{-1}-2k_2 .
\label{eq:B3}
\end{equation}
Unlike $B^{(1)}$, it is not a translation bound. The term $q_1$ is not convex when $\tau > k_1$, and the minimization in $S_1$ uses the restriction \eqref{eq:psd}.
The field $E_1$ is rank one, while the Lagrangian of $S_2$ is degenerate and only the trace of $E_2$ is fixed,
\[
\Lambda_1=\Bigl\{ (r/(2k_1),0,0)\Bigr\}, \quad
\Lambda_2=\Bigl\{\lambda\ge0:\ \sum_j\lambda_j=\frac{3}{a_*}\Bigr\},\qquad a_*=\frac{k_2}{k_1}m_1+m_2 ,
\]
and $a_*$ is the constant of \eqref{eq:sets3}. As $m_1\to0$ the bound \eqref{eq:B3} tends to
$k_2(3-2m_2)/m_2$, the two-material Hashin--Shtrikman value for an ideal conductor in a matrix of
$S_2$, as expected. At the upper end $\xi_4$ the multiplier $\tau$ reaches $k_2$, where
$1/\kappa=6k_2$. As at $\xi_3$, $r$ is the largest positive root of the cubic in \eqref{eq:xi4}, and the first equation of \eqref{eq:pair} gives $\xi_4$:
\begin{equation}
m_2r^{3}-48k_2^{2}r+72k_2^{3}=0,\qquad \xi_4=2k_1\Bigl(\frac{3}{r}-\frac{m_2}{2k_2}\Bigr).
\label{eq:xi4}
\end{equation}
With $r=6k_2/a_*$, which is the value of $r$ in this region, the cubic becomes
$a_*^3-4a_*^2+3m_2=0$, and the region is
\begin{equation}
a_*^2(4-a_*)\le3m_2 ,
\label{eq:cond3d}
\end{equation}
the spatial form of \eqref{eq:cond3}.
\begin{remark}
Briane and Nesi [2004, Th.~4.4] proved the bound $B^{(3)}$, as
the limit $k_3\to\infty$ of their bound with the translation $k_2$, for the quasiorthogonal composites (Section~\ref{sec:intro}), a class larger than that of \eqref{eq:psd} that contains the coaxial laminates [Briane and Nesi 2004, Th.~2.15]. Their coaxial laminates of rank 3 attain it isotropically only for
$m_1\le(k_1/k_2)m_2(1-m_2)/(2+m_2)$ [Briane and Nesi 2004, Cor.~4.6], below $\xi_4$; the
structures of Section~\ref{sec:structures3d} attain it on the whole region.
\end{remark}

\begin{figure}[htbp]
\centering
\IfFileExists{nov18b_fields3d.png}{\includegraphics[width=0.9\linewidth]{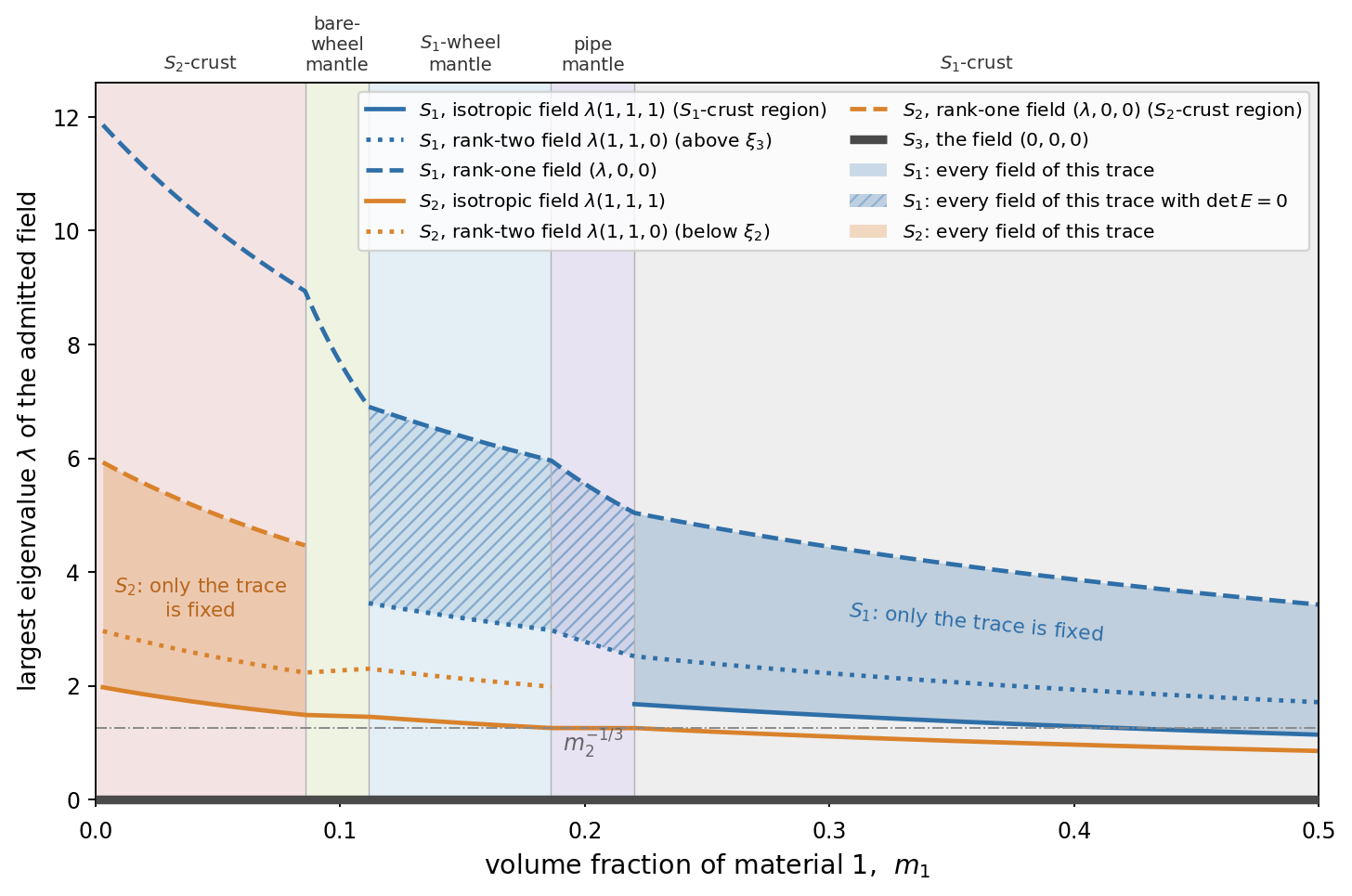}}{%
\framebox[0.8\linewidth][c]{\rule{0pt}{4cm}\textit{[figure file not found]}}}
\caption{The fields that the bound admits in $S_1$ (blue), $S_2$ (orange) and $S_3$ (dark), drawn by
their largest eigenvalue $\lambda$: solid for the isotropic, dotted for the rank-two and dashed for
the rank-one fields; the shading spans the admitted families. $k_1=1$, $k_2=2$, $m_2=\tfrac12$.}
\label{fig:fields3d}
\end{figure}

The optimal fields are shown in Figure~\ref{fig:fields3d}. The field of $S_1$ has the trace
$r/(2k_1)$ in every region. In the $S_1$-crust region every field of $S_1$ of this trace is
admitted, from the isotropic field to the rank-one field, and the shading spans the largest
eigenvalues of these fields; in the pipe mantle and $S_1$-wheel mantle regions only the fields with
$\det E_1=0$ are admitted, from $\lambda(1,1,0)$ to $(\lambda,0,0)$ (hatched); below $\xi_3$ the
field of $S_1$ is rank one. The field of $S_2$ is isotropic in the $S_1$-crust and pipe mantle regions, with $\lambda=m_2^{-1/3}$ (dash-dot line) independently of $m_1$ in the second; between $\xi_2$ and $\xi_4$ it takes two values, the isotropic field and the rank-two field
\eqref{eq:twoatom}; in the $S_2$-crust region every field of $S_2$ of the trace $3/a_*$ is admitted.
Each curve is continuous across the boundaries it crosses.

\begin{figure}[htbp]
\centering
\IfFileExists{nov18_supportsets3d.png}{\includegraphics[width=\linewidth]{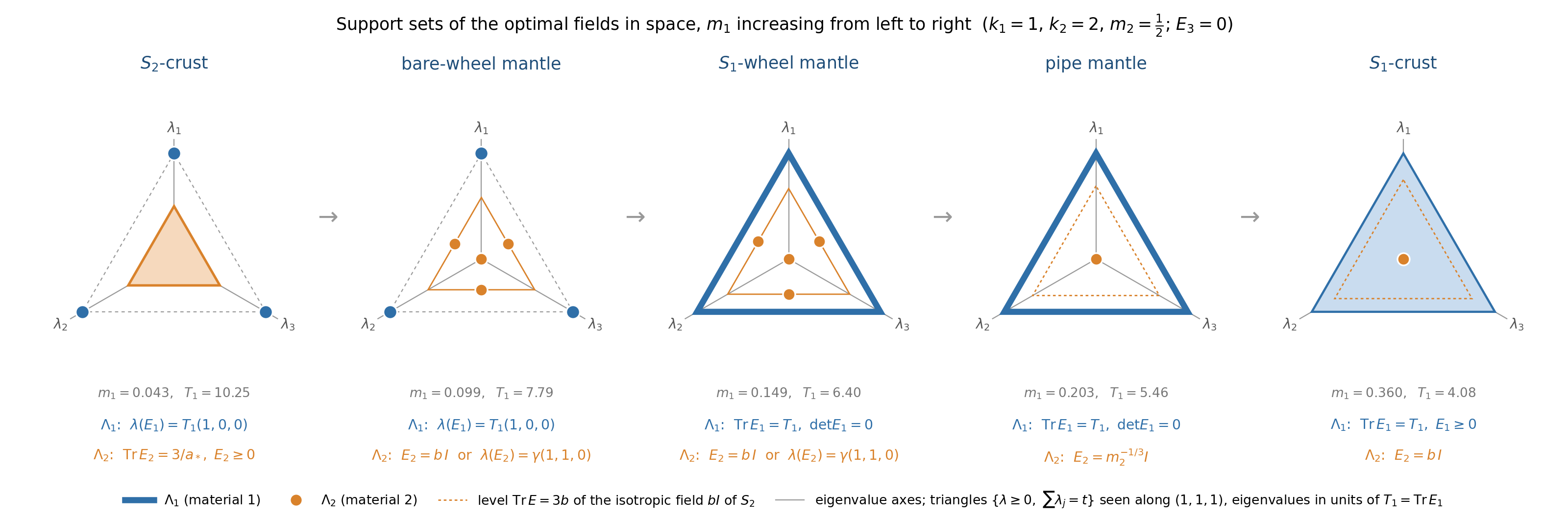}}{%
\framebox[0.8\linewidth][c]{\rule{0pt}{4cm}\textit{[figure file not found]}}}
\caption{The support sets of the three-dimensional bound, one value of $m_1$ in each region;
$k_1=1$, $k_2=2$, $m_2=\tfrac12$.}
\label{fig:sets3d}
\end{figure}

Figure~\ref{fig:sets3d} illustrates the support sets in each of the five regions. Each panel shows
the octant of eigenvalues $\lambda\ge0$ seen along the diagonal $(1,1,1)$, so that a plane of
constant trace appears as a triangle with its vertices on the axes; the eigenvalues are in units of
$T_1=\Tr E_1=r/(2k_1)$. The support set $\Lambda_1$ of $S_1$ (blue) is the whole triangle
$\Tr E_1=T_1$ in the $S_1$-crust region, its edges ($\det E_1=0$) in the pipe mantle and $S_1$-wheel mantle regions, and its vertices (rank-one fields) below $\xi_3$. The support set $\Lambda_2$ of
$S_2$ (orange) is the isotropic point $bI$ (dashed in the two right panels: the triangle of its trace $3b$), together with
the three rank-two points $\gamma(1,1,0)$ between $\xi_2$ and $\xi_4$ (on the solid triangle of
the trace $2\gamma$), and the whole triangle $\Tr E_2=3/a_*$ in the $S_2$-crust region. Everywhere
$E_3=0$. Compare the planar Figure~\ref{fig:sets2d}.

The five regions are the intervals of $m_1$ on which the same restrictions are active.
At their boundaries the type of the minimizing field
changes in one of the materials (Figures~\ref{fig:fields3d} and~\ref{fig:sets3d}). Two mechanisms are
responsible. The first is the planar one, the multiplier $\tau$ against the conductivities: as
$m_1$ decreases, $\tau$ moves from $k_1$ to $k_2$. It leaves $k_1$ at $\xi_3$, where $S_1$ passes from a field with only the trace and $\det E_1=0$ fixed to a strictly rank-one field. It
reaches $k_2$ at $\xi_4$, where the Lagrangian of $S_2$ becomes degenerate.

The second mechanism is the integral restriction on the determinant. Above $\xi_1$ the
Hashin--Shtrikman fields satisfy it, because $S_1$ can lower its own determinant to zero. Below
$\xi_1$ they cannot; then $\mu$ becomes positive, $S_1$ is driven to $\det E_1=0$, and the restriction rests on $S_2$ alone, which fixes its isotropic field at $E_2=I/a$. At $\xi_2$
the rank-two field $(\gamma,\gamma,0)$ ties with the isotropic one, and between $\xi_2$ and $\xi_4$ the material $S_2$ carries both, in the proportion set by this restriction. At $\xi_4$ the
multiplier $\mu$ returns to zero, together with $\tau=k_2$.

All four transitions are thresholds of $r$, which decreases monotonically with $m_1$. Indeed, the
bound is a maximum of functions affine in $m_1$, hence convex in $m_1$, and its derivative
$dk_*/dm_1=g_1/3=-r^2/(12k_1)$ does not decrease. The threshold values are $r=2(k_2+2k_1)/a$ at
$\xi_1$, $r=(2C+\sqrt{2Ck_2})/a$ at $\xi_2$, and the largest positive roots of the cubics
\eqref{eq:xi3} and \eqref{eq:xi4} at $\xi_3$ and $\xi_4$.

\begin{remark}[The order of the thresholds]\label{rem:order}
The threshold values of $r$ are ordered so that $\xi_1>\xi_2>\xi_3>\xi_4$ for all
$k_1<k_2$ and $0<m_2<1$. First, $\xi_2<\xi_1$ by \eqref{eq:xi2}, since
$\sqrt{2k_2(k_1+k_2)}>2k_1$. Next, $2C-A=k_2$ gives $1/\kappa=2C+\sqrt{2Ck_2}$ at every $\tau$, so
the value at $\xi_2$ is $r=1/(\kappa a)$, where $p=1$ in \eqref{eq:twoatom}. Let
\[
f_\tau(r)=m_2r^3-12C^2r+\frac{12C^2}{\kappa}-\frac{1}{\kappa^3}
\]
be the second equation of \eqref{eq:pair} multiplied by $2C^2r$; at $\tau=k_1$ it is the cubic of
\eqref{eq:xi3}, and at $\tau=k_2$, where $1/\kappa=6k_2$, the cubic of \eqref{eq:xi4}. With
$t=\sqrt{2k_2/C}\in[1,\sqrt2)$, its free term $C^3(2+t)(8-4t-t^2)$ is positive, while
$f_\tau\bigl(1/(\kappa a)\bigr)=12C^2(1-1/a)/\kappa<0$. Hence $f_\tau$ has two positive roots, and
only the larger one gives $p<1$; at $\tau=k_1$ it exceeds $1/(\kappa a)$, so $\xi_3<\xi_2$.
Finally, the derivative of the free term in $\tau$,
$C^2\bigl[(2+t/2)(8-4t-t^2)+(2+t)(16-6t-t^2)\bigr]\le C^2\bigl[7.5+9(2+t)\bigr]$, is less than
$24C^2(2+t)=24C/\kappa$. Thus $\partial f_\tau/\partial\tau<24C/\kappa-24Cr<0$ for $r>1/\kappa$,
and the larger root increases with $\tau$. Its value at $\tau=k_2$ exceeds that at $\tau=k_1$,
so $\xi_4<\xi_3$.
\end{remark}

The bounds of the two outer regions have the form of the Hashin--Shtrikman bound,
\begin{equation}
\frac{1}{B^{(1)}+2k_1}=\frac{m_1}{3k_1}+\frac{m_2}{k_2+2k_1}+0\cdot m_3,\qquad
\frac{1}{B^{(3)}+2k_2}=\frac{m_1}{3k_1}+\frac{m_2}{3k_2}+0\cdot m_3 ,
\label{eq:hsform3d}
\end{equation}
affine in the fractions. A sphere of conductivity $k_c$ coated by a shell of $S_j$ that takes the fraction $\beta$ of its volume is neutral in the medium of conductivity $k$ when \eqref{eq:coat} holds with $d=3$. The weight $1/(3k_j)$ of the coat equals the weight of $S_1$ in the first bound and of $S_2$ in the
second, exactly as in the plane. Thus the enveloping principle holds in space as in the plane:  if a structure
attains the bound $B^{(1)}$ at $\hat m$, the assemblage of its spheres coated by $S_1$ attains it at
every point of the ray toward pure $S_1$ that lies in the $S_1$-crust region; the same holds for a coat of $S_2$ in the
$S_2$-crust region.

\begin{figure}[htbp]
\begin{center}
\IfFileExists{nov18b_regions3d.png}{\includegraphics[width=0.95\linewidth]{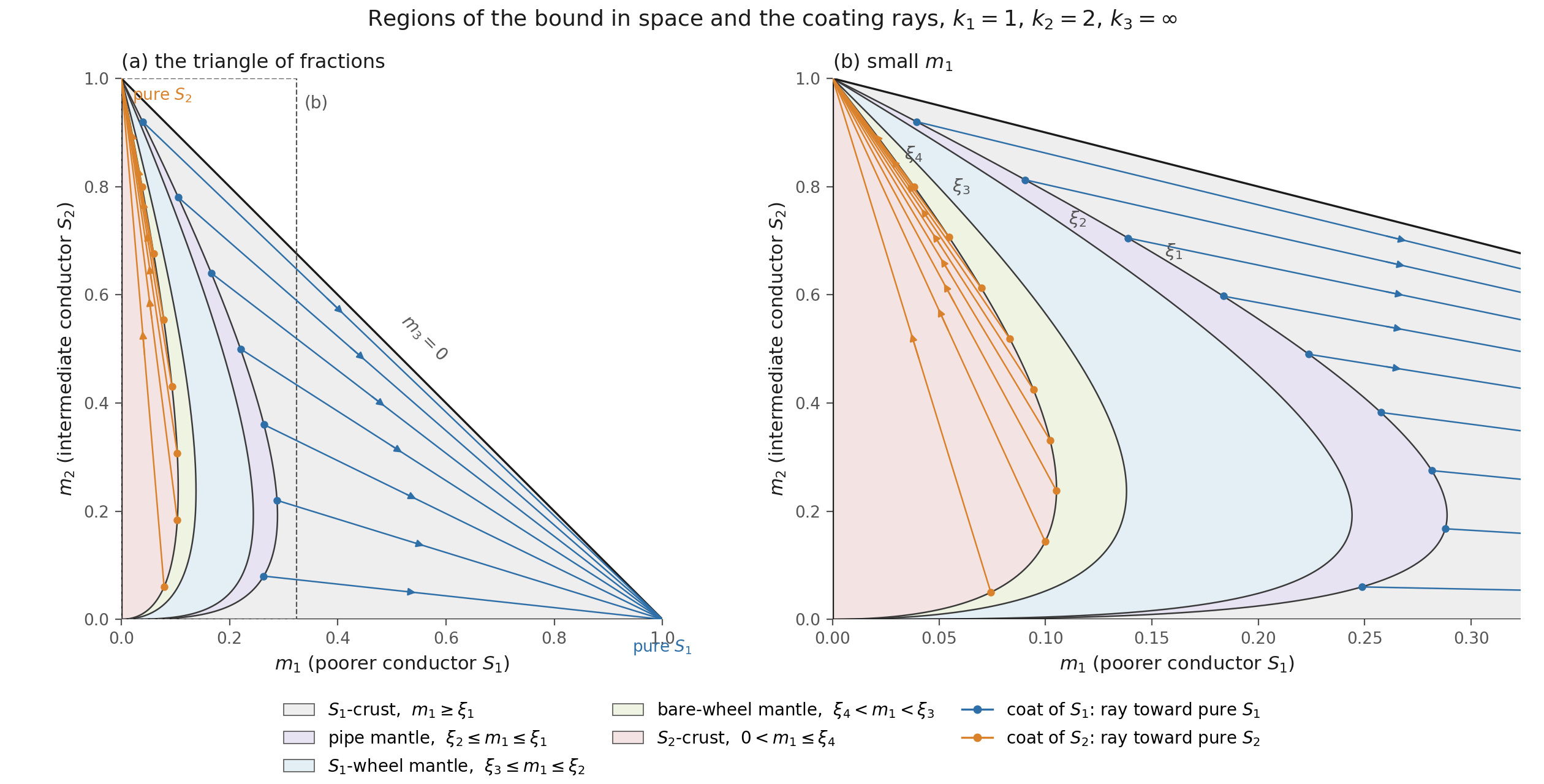}}{%
\framebox[0.8\linewidth][c]{\rule{0pt}{4cm}\textit{[figure file not found]}}}
\caption{The five regions of the three-dimensional bound and the coating rays in the plane of the
fractions $(m_1,m_2)$, $k_1=1$, $k_2=2$: (a) the whole triangle, (b) the strip of small $m_1$.}
\label{fig:regions3d}
\end{center}
\end{figure}

As in the plane, 
every point of the $S_1$-crust region other
than the vertex $m_1=1$ lies on exactly one segment from the curve $m_1=\xi_1(m_2)$ to this
vertex. Likewise, every point of the $S_2$-crust region other than the vertex $m_2=1$ lies on exactly one
segment from the curve $m_1=\xi_4(m_2)$ to this vertex (Figure~\ref{fig:regions3d}). Therefore
\emph{optimal structures need to be found only in the three intermediate regions}, including their
end points, and the two outer regions follow.

In Figure~\ref{fig:regions3d}, as in Figure~\ref{fig:regimes2d}, the blue rays from the curve $m_1=\xi_1$ fill the $S_1$-crust region, and the orange rays from the curve $m_1=\xi_4$ fill the $S_2$-crust region.

\subsection{Optimal structures}\label{sec:structures3d}

In this section we construct, region by region, structures that attain the bound.
Each atom of the structure consists of a spherical core of $S_2$ surrounded by a spherical layer,
the mantle. The mantle is an optimal planar structure, extended along the radius and wrapped around
the core: the planar Hashin--Shtrikman assemblage of discs of $S_3$ coated by $S_1$, or one of the
wheels of Section~\ref{sec:structures2d}. In the outer regions, with extreme $m_1$, an additional
layer, the crust of $S_1$ or $S_2$, is added. Notice that the corresponding planar assemblage is
similar: it contains the circular core of $S_2$, surrounded by the optimal one-dimensional
structure (a laminate) coiled around it, plus the outer crust of $S_1$ or $S_2$ in the regimes with extreme $m_1$.

\begin{figure}[htbp]
\begin{center}
\IfFileExists{nov18b_summary3d.png}{\includegraphics[width=\linewidth]{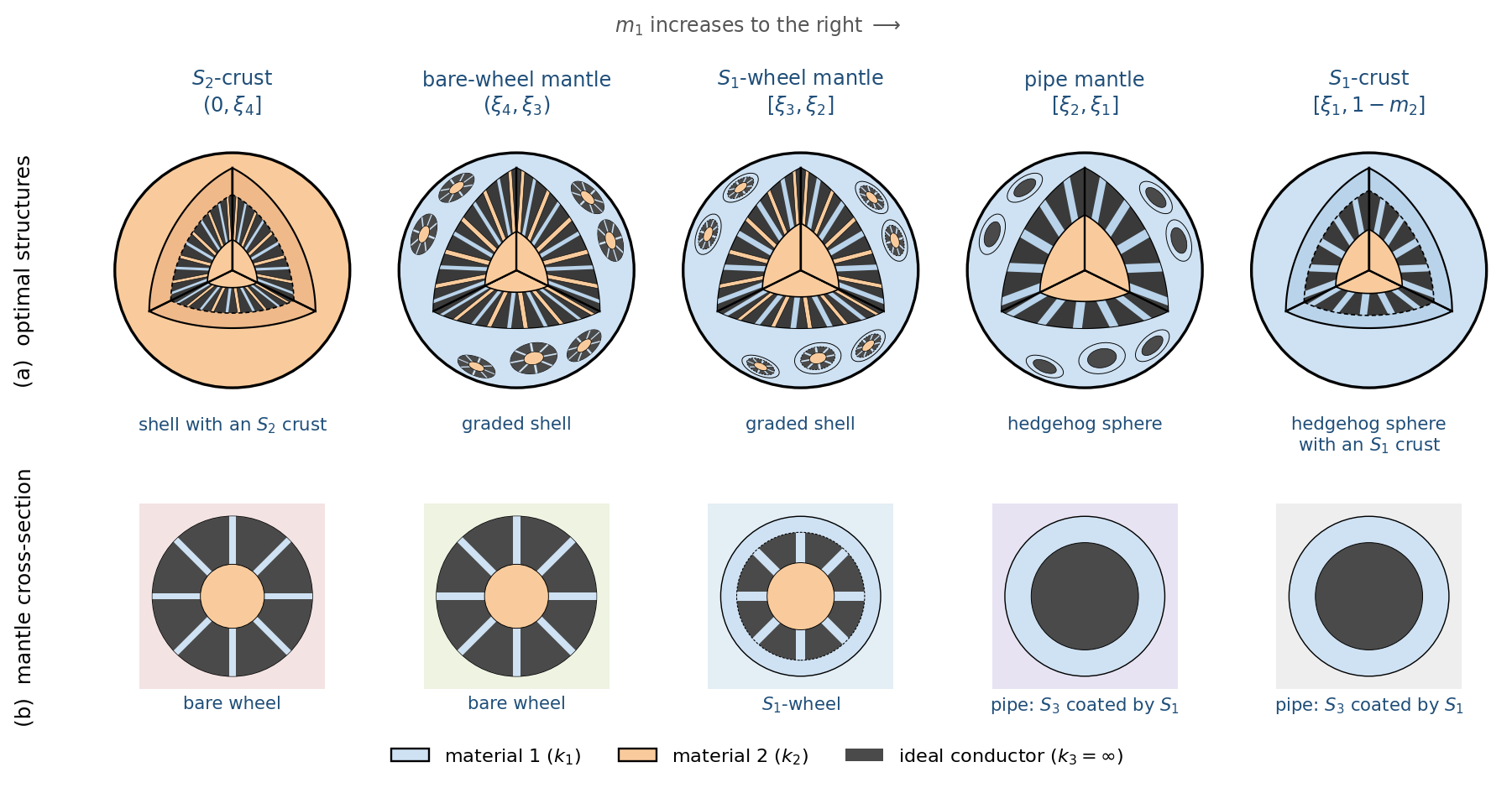}}{%
\framebox[0.8\linewidth][c]{\rule{0pt}{4cm}\textit{[figure file not found]}}}
\caption{The optimal structures in space, $m_1$ increasing to the right: (a) the atoms with the
near octant removed; (b) the cross-section of the mantle.}
\label{fig:summary3d}
\end{center}
\end{figure}

Figure~\ref{fig:summary3d} shows these atoms, $m_1$ increasing to the right (from the left: the
$S_2$-crust, bare-wheel mantle, $S_1$-wheel mantle, pipe mantle and $S_1$-crust regions). Panel (a)
shows the structures as spheres with the near octant removed: the graded shell of the bare-wheel
mantle region with a crust of $S_2$; the graded shells whose cross-section is the planar bare wheel
and the planar $S_1$-wheel; the hedgehog sphere; and the hedgehog sphere with a crust of $S_1$. On the cut faces the layers
of $S_1$ and the rods of $S_2$ keep a constant width, and the conductor between them widens
outward. Panel (b) shows the mantle in cross-section, the planar pattern carried by every
concentric sphere of the mantle: the bare wheel (in the bare-wheel mantle region, and under the
crust in the $S_2$-crust region), whose disc of $S_2$ is the section of a radial rod; the
$S_1$-wheel (in the $S_1$-wheel mantle region); and the disc of $S_3$ coated by $S_1$, the section of
a radial pipe (no $S_2$ in the mantle: in the pipe mantle region, and under the crust in the
$S_1$-crust region). These structures explain the names of the regions: in the three intermediate
regions the mantle consists of radial pipes, of extruded $S_1$-wheels or of extruded bare wheels,
and in the two outer regions a crust of $S_1$ or of $S_2$ is added.

\paragraph{Fields and geometric elements.} Let us begin with
the fields that the bound admits and with the geometric elements that realize them. The support fields of \eqref{eq:g} are of three
types: isotropic, $(b,b,b)$; rank two, $(\gamma,\gamma,0)$, for $S_2$; and rank one,
$(\lambda,0,0)$. Where a Lagrangian is degenerate, the minimizers form whole families of one trace.
If $\tau=k_1$ and $\mu>0$, these are the fields of $S_1$ with $\det E_1=0$ (the edges of the
triangle). If $\mu=0$, they are all the fields of the prescribed trace ($S_1$ in the
$S_1$-crust region, $S_2$ in the $S_2$-crust region). There are many optimal structures; we use
radial ones for illustration, because they are isotropic by construction. In them the isotropic,
rank-two and rank-one fields correspond to three geometric elements.
\begin{itemize}
\item An \emph{isotropic} field is carried by a region that the ideal conductor does not cut. In the
radial structures this is the central sphere of $S_2$, of radius $a_h$. The conductor joins its
surface to the outer surface along the radii and fixes its field at $I/a_h$.
\item A field with one zero eigenvalue is carried by a \emph{perforated sheet}, a sheet of $S_1$ or $S_2$ penetrated by thin radial cylinders of $S_3$. 
The conductor $S_3$ keeps the potential constant along
the radius, so the radial eigenvalue vanishes, and the field in the transverse plane is the planar
field of the cross-section.
\item A \emph{rank-one} field is carried by a \emph{sheet} neighboring $S_3$, a thin radial layer
between two $S_3$  faces with the normal $n$. Across each face the field jumps from $E_3=0$ by a
rank-one matrix, so in the sheet $E=\lambda\,n\otimes n$, as in Section~\ref{sec:structures2d}.
\end{itemize}
Notice that the three elements form a path in the rank of the field. Starting from $E_3=0$ in the
conductor, one rank-one jump gives the field of a sheet, two rank-one terms give the field of a perforated sheet,
and three give the isotropic field of the sphere. 

%

\paragraph{The common topology.} In the three intermediate regions the optimal structures that we
use for illustration are radial assemblages of two layers. (i) The core is the sphere of $S_2$ of radius $a_h$, with the isotropic field $I/a_h$. (ii) The mantle is the shell $a_h<\rho<1$. In it the conductor $S_3$ runs along
the radii from the core to the outer surface. The cross-section of the mantle at every radius is an
optimal planar structure at the local fractions, extruded along the radius
(Figure~\ref{fig:summary3d}(b)). In the two outer regions a crust from either $S_1$ or $S_2$ is added
(Section~\ref{sec:envelopes3d}). We write $\rho$ for the radius, keeping $r$ for the multiplier, and
$m_{1L}(\rho)$, $m_{2L}(\rho)$ for the local fractions. The volume restrictions are
\begin{equation}
3\int_{a_h}^{1}m_{1L}(\rho)\,\rho^2\,d\rho=m_1,\qquad
a_h^3+3\int_{a_h}^{1}m_{2L}(\rho)\,\rho^2\,d\rho=m_2 .
\label{eq:volumes}
\end{equation}

\paragraph{The radial quadrature.} Let us compute the conductivity $s$ of a radial assemblage. Grow
its outer radius from $\rho$ to $\rho(1+\delta)$ by a layer that conducts infinitely along the radius
and with $k_t$ in the two tangential directions. The field in the layer is purely tangential, and
after rescaling the new conductivity is $(s+2k_t\delta)/(1+\delta)$. So $ds=\delta(2k_t-s)$, while
the volume $v=\rho^3$ grows by $dv/v=3\delta$. Hence
\begin{equation}
\frac{ds}{dv}=\frac{2k_t-s}{3v},\qquad\text{that is,}\qquad \frac{d(\rho s)}{d\rho}=2k_t.
\label{eq:ds}
\end{equation}
 For a graded sphere with the radial conductivity $\lambda_n$
and the tangential conductivity $\lambda_t$ (here
$\lambda_n$ and $\lambda_t$ are conductivities, not eigenvalues of the field), Tartar
found the Riccati equation $d(\rho s)/d\rho+s^2/\lambda_n=2\lambda_t$ [Milton 2002, \S7.4]. Equation
\eqref{eq:ds} is its case $\lambda_n=\infty$, $\lambda_t=k_t$. Integrating
\eqref{eq:ds} from the central sphere, we obtain
\begin{equation}
s(1)=k_2a_h+2\int_{a_h}^{1}k_t(\rho)\,d\rho .
\label{eq:quadrature}
\end{equation}
Physically, the central sphere carries the field $I/a_h$ and contributes $k_2a_h$. The radial
direction of the mantle carries no energy, so each layer contributes only through its transverse
conductivity $k_t(\rho)$, the conductivity of the planar cross-section.

Minimizing \eqref{eq:quadrature} over the profiles under the volume restrictions \eqref{eq:volumes}
is a pointwise isoperimetric problem. Its Euler--Lagrange conditions, $\partial k_t/\partial
m_{iL}\propto\rho^2$ with one constant per restriction, give the optimal profiles. Equivalently, as
in the plane, the profiles keep the fields at the values that the bound prescribes. Strictly
speaking, the radial structures are limits of iterated homogenization with a cross-section that
changes with the radius; such limits are realized by sequences of microstructures [Avellaneda 1987].

We now construct the structures, starting from the pipe mantle region, the second panel from the right of Figure~\ref{fig:summary3d}.

\subsubsection{The pipe mantle: the hedgehog sphere}\label{sec:spiky}

In the pipe mantle region all of $S_2$ is in the central sphere, whose radius $a=m_2^{1/3}$ does not
depend on $m_1$, and the mantle contains only $S_1$ and $S_3$. The structure is a solid sphere of
$S_2$ surrounded by a mantle of radial ideal-conductor spikes in $S_1$. We call it the
\emph{hedgehog sphere}. Its mantle is a radially oriented uniaxial shell, of the type introduced by
Schulgasser [1983] (see Milton [2002, \S7.5]). With a crust of $S_1$, the hedgehog sphere at its
threshold $\xi_1$ is the ``hairy sphere'' of Albin and Cherkaev [2006] at $k_3=\infty$
(Section~\ref{sec:envelopes3d}). New here is
the optimality of the hedgehog sphere in the pipe mantle region, where the bound exceeds the
Hashin--Shtrikman bound. Let $\nu(\rho)=m_{1L}(\rho)$ be the local fraction of $S_1$, with
$3\int_a^1\nu\rho^2\,d\rho=m_1$.

\paragraph{The grading is fixed by the field.} At the radius $\rho$ the average tangential field is
$1/\rho$ in each of the two tangential directions. The field in $S_3$ is zero, so the mean field in
$S_1$ is $1/(\rho\nu)$ in each tangential direction, and the radial one is zero. Therefore
$\det E_1=0$ and $\Tr E_1=2/(\rho\nu)$. The bound prescribes a constant trace \eqref{eq:fields-sb},
so $\nu\propto 1/\rho$, and the volume restriction of $S_1$ gives
\begin{equation}
\nu(\rho)=\frac{C_w}{\rho},\qquad C_w=\frac{2m_1}{3(1-a^2)},\qquad
\Tr E_1=\frac{2}{C_w}=\frac{3(1-a^2)}{m_1}.
\label{eq:nu}
\end{equation}

The same profile follows from optimization. The cross-section is a planar composite of $S_3$
inclusions in $S_1$, whose least conductivity is $k_t=2k_1/\nu-k_1$ (the limit $m_2=0$ of
\eqref{eq:k01}). Minimizing $k_2a+2\int_a^1(2k_1/\nu-k_1)\,d\rho$ at the fixed volume, we obtain
$2k_1/\nu^2\propto\rho^2$. An optimal cross-section is the assemblage of discs of $S_3$ coated by
$S_1$ (Figure~\ref{fig:summary3d}(b)); in space these are radial conductor spikes, each in a sleeve of $S_1$, the pipes that give the region its name. Physically, the conductor spikes carry the potential of the outer surface to the core without loss,
and the tangential current has to cross the sleeves of $S_1$ between neighboring spikes.

\paragraph{Bound matches structure.} Let us show that the hedgehog sphere realizes the bound.
Substituting \eqref{eq:nu} into \eqref{eq:quadrature}, we obtain
\[
k_S=k_2a+\frac{2k_1}{C_w}\bigl(1-a^2\bigr)-2k_1(1-a)=k_2a+\frac{3k_1(1-a^2)^2}{m_1}-2k_1(1-a)=B_2^{*},
\]
which is \eqref{eq:B2}. The fields also coincide with those of the bound: the field in the sphere is
$I/a$, and the trace of $E_1$ in \eqref{eq:nu} is $r/(2k_1)$ at $r=6k_1(1-a^2)/m_1$. This trace holds
at every point of $S_1$, not only on average. Indeed, in the coat of a coated disc the field is
$c_1I$ plus a traceless harmonic part, so its trace is constant. Its eigenvalues are nonnegative,
and the radial one is zero.

\subsubsection{The \texorpdfstring{$S_1$}{S1}-wheel mantle}

Below $\xi_2$ the bound prescribes two values of the field in $S_2$.
Accordingly, the central sphere shrinks to a radius $a_h<m_2^{1/3}$, and the rest of
$S_2$ moves into the mantle as radial rods, where its field is $(\gamma,\gamma,0)$. The mantle is now
a graded three-material shell. Its cross-section at every radius is an optimal planar
three-material structure at the local fractions. As $m_1$ decreases, the cross-section passes through the planar $S_1$-wheel and bare-wheel regimes and, at $\xi_4$, reaches the $S_2$-wheel regime. The thresholds $\xi_3$ and $\xi_4$ are the fractions at
which it changes the planar regime at every radius simultaneously.

The field in the sphere is $b\,I$ with $b=1/a_h=\kappa r$, and the field in the rods is
$(\gamma,\gamma,0)$ with $\gamma=r/(2C)$, the same at every radius. Notice that the two values have equal translated energy but not equal traces, since $\det E_2$ is zero in the rods. Material $S_1$ is only in the mantle, in the
sheets around the rods of $S_2$ and in the cylindrical sleeve of $S_1$ around them (the extruded crust of the planar $S_1$-wheel). Its determinant is zero everywhere, and its trace
is constant.

The cross-section lies on the planar branch of the $S_1$-wheel regime, $k_t=\Sigma(\rho)^{-1}-k_1$, where
\[
\Sigma(\rho)=\frac{m_{1L}(\rho)}{2k_1}+\frac{m_{2L}(\rho)}{k_1+k_2}.
\]
Both Euler--Lagrange conditions say that $\Sigma^{-2}\propto\rho^2$, so the optimal mantle obeys the
local law $\Sigma(\rho)=K/\rho$. Only the combination $\Sigma$ is fixed. The split of $\Sigma(\rho)$
between $S_1$ and $S_2$ is free, subject to the ceiling $\sqrt{m_{2L}(\rho)}\le(k_1+k_2)K/\rho$ that
keeps the cross-section in the planar $S_1$-wheel regime. The volume restrictions \eqref{eq:volumes} give
$K=2\bar\Sigma(a_h)/(3(1-a_h^2))$, and \eqref{eq:quadrature} gives the conductivity
\begin{equation}
k_S(a_h)=k_2a_h-2k_1(1-a_h)+\frac{3(1-a_h^2)^2}{2\bar\Sigma(a_h)},\qquad
\bar\Sigma(a_h)=\frac{m_1}{2k_1}+\frac{m_2-a_h^3}{k_1+k_2}.
\label{eq:kS1}
\end{equation}
The only free parameter left is the radius $a_h$ of the central sphere. The conductivity
\eqref{eq:kS1} is stationary in $a_h$ where
\begin{equation}
(2\omega-1)\,a_h^3+3a_h=(1+\omega)\Bigl(\frac{m_1(k_1+k_2)}{k_1}+2m_2\Bigr),\qquad
\omega=\sqrt{\frac{k_2}{2(k_1+k_2)}},
\label{eq:ahcubic}
\end{equation}
a cubic with a single root in $(0,1)$. At $m_1=\xi_2$ the root is $a_h=m_2^{1/3}$, and the
structure becomes the hedgehog sphere.

\paragraph{Bound matches structure.} At the root of \eqref{eq:ahcubic} the value \eqref{eq:kS1}
equals the bound \eqref{eq:wsbound} identically. Indeed, take the trial multiplier
$r=3(1-a_h^2)/\bar\Sigma(a_h)$. The difference of the two closed forms is then a rational function of
$a_h$ whose numerator contains the cubic \eqref{eq:ahcubic} twice, while the cubic \eqref{eq:cubic}
of the bound is divisible by it. 
The region ends at $\xi_3$, where the ceiling on $m_{2L}$ in the mantle is reached at every
radius at once. At saturation the profile is exactly the profile of the next region, with its local
planar multiplier equal to $k_1$.

\subsubsection{The bare-wheel mantle}\label{sec:bare}

As $m_1$ decreases further, the mantle loses the sleeve of $S_1$ around the rods, and its
cross-section becomes the planar bare wheel. Material $S_2$ is in two locations, the central
sphere and the radial rods, which are the extruded discs of the planar atom. As in the previous
region, the translated energy of $E_2$ is the same in both. The fields in $S_2$ are as in the
previous region, and the fields in $S_1$ are rank one with a constant magnitude.

The cross-section lies on the planar branch of the bare-wheel regime, $k_t=k_2+2k_1(1-\sqrt{m_{2L}})^2/m_{1L}$.
Now the two Euler--Lagrange conditions determine both local fractions,
\begin{equation}
\sqrt{m_{2L}(\rho)}=\frac{u}{\rho},\qquad
m_{1L}(\rho)=\frac{C_1}{\rho}\Bigl(1-\frac{u}{\rho}\Bigr),\qquad
k_t(\rho)=k_2+\frac{2k_1}{C_1}(\rho-u),
\label{eq:profile2}
\end{equation}
with two constants $u$ and $C_1$ fixed by the volume restrictions. Notice that the rods of $S_2$ have
a constant cross-section ($m_{2L}\propto1/\rho^2$). The volume restrictions \eqref{eq:volumes} give
\[
3u^2(1-a_h)=m_2-a_h^3,\qquad C_1=\frac{2m_1}{3W_1},\qquad W_1=(1-a_h)(1+a_h-2u),
\]
and \eqref{eq:quadrature} gives
\begin{equation}
k_S(a_h,u)=k_2(2-a_h)+\frac{3k_1W_1^2}{m_1}.
\label{eq:kS2}
\end{equation}
Its stationarity under the volume restriction of $S_2$ reads
\begin{equation}
k_2m_1u=6k_1(1-a_h)(a_h-u)^2(1+a_h-2u).
\label{eq:ahu}
\end{equation}

\paragraph{Bound matches structure.} The local planar multiplier of the cross-section,
\[
k_2+\tau=\frac{2k_1\bigl(\sqrt{m_{2L}}-m_{2L}\bigr)}{m_{1L}}=\frac{2k_1u}{C_1},
\]
does not depend on $\rho$. Thus one interior value of $\tau$ serves the whole mantle, and the
cross-section is in the planar bare-wheel regime at every radius exactly when $\tau\in[k_1,k_2]$. The map
\[
r=\frac{6k_1W_1}{m_1},\qquad k_2+\tau=\frac{r\,u}{2}
\]
carries the $2\times2$ system of the structure (the volume restriction of $S_2$ and \eqref{eq:ahu})
onto the system \eqref{eq:pair} of the bound. The identity of \eqref{eq:kS2} and \eqref{eq:wsbound} is proved for all parameters on the branch
of the region, $0<u\le2a_h/3$ and $\tau\ge k_1$ (the first condition is $\tau\le k_2$). Indeed, the
system of the structure contains the fractions linearly. Solving it for them, we obtain
\[
m_2=a_h^3+3u^2(1-a_h),\qquad m_1=\frac{6k_1(1-a_h)(a_h-u)^2(1+a_h-2u)}{k_2u},
\]
after which the difference between \eqref{eq:kS2} and the bound \eqref{eq:wsbound} cancels
identically in $(a_h,u,k_1,k_2)$. 
When $\tau$ reaches $k_2$, the cross-section reaches the boundary of the planar $S_2$-wheel regime at every
radius at once; this is the threshold $\xi_4$.

\subsubsection{The outer regions: the crust}\label{sec:envelopes3d}

By the enveloping principle, the structures of the two outer regions are the structures at the ends
of the intermediate range, with a crust.

In the $S_1$-crust region the structure is the hedgehog sphere at its own threshold $\xi_1(\hat m_2)$,
coated by $S_1$. The central sphere holds all of $S_2$ and has the radius $(1-\beta)^{1/3}\hat m_2^{1/3}=m_2^{1/3}$,
independently of $\beta$; only the crust changes. When $m_3=0$, the structure is the two-material
Hashin--Shtrikman coated sphere of $S_2$ in $S_1$. As $m_2\to0$, the hedgehog sphere under the crust
becomes pure conductor, and the structure becomes the coated sphere of $S_3$ in $S_1$, with the
conductivity $k_1(3-2m_1)/m_1$.

In the $S_2$-crust region, the graded shell of the bare-wheel mantle region at its threshold
$\xi_4(\hat m_2)$ is coated by $S_2$. Let us parametrize the threshold curve by $a_*\in(0,1)$. On it
the fractions are $\hat m_2=a_*^2(4-a_*)/3$ and $\hat m_1=\xi_4=(k_1/k_2)\,a_*(1-a_*)(3-a_*)/3$,
and the ray parameter $\hat m_1/(1-\hat m_2)=(k_1/k_2)\,a_*(3-a_*)/(3+3a_*-a_*^2)$ increases with
$a_*$. The structure under the crust at a given point is the one where its ray meets the curve. As
$m_1\to0$, this structure becomes pure conductor, and the whole structure becomes the coated sphere
of $S_3$ in $S_2$, with the conductivity $k_2(3-2m_2)/m_2$.

The bounds $B^{(1)}$ and $B^{(3)}$ are values of
the bound functional at every fraction, and the coating identity is exact. The structures under the
crust attain the bound at their own fractions. In the $S_1$-crust region the structure with a crust of $S_1$ attains $B^{(1)}$, and in the $S_2$-crust region the structure with a crust of $S_2$ attains $B^{(3)}$.


\section{Discussion}\label{sec:discussion}

In this work, we found the lower bounds on the effective properties of the composite made of
three conductors, one of them an ideal conductor. The bounds have three regimes in the plane and five regions in space; in each of them we found optimal structures and proved analytically that their conductivity equals the bound. New are the bound in space between the
$S_1$-crust region and the $S_2$-crust region, and isotropic optimal structures on the whole range
below the $S_1$-crust region. 
There are many questions that remain open.

The main open question is Conjecture (H). Cherkaev [2024] noted that a systematic derivation of
the pointwise restrictions on the fields is still missing. The bound is proved for every structure when
$m_1\ge\xi_1$, and below $\xi_1$ for the structures whose fields satisfy \eqref{eq:psd}, among
them all the structures constructed here. That it holds for all structures follows from (H),
which at present is open.

In this work, we assumed that $S_3$ is an ideal conductor. For a finite $k_3$ the field in $S_3$ is not zero,
the argument for \eqref{eq:psd} based on the compatibility with $E_3=0$ is lost, and the bounds based on the restriction \eqref{eq:psd} are not exact [Cherkaev 2009, \S4.4]. The open question is what type of restrictions is
needed to find the exact bounds in the case of finite conductivities.

The method extends formally to more materials and to anisotropic composites; we leave
these extensions to future work. Our first attempts show that the bounds can be obtained, but
there is a gap between the bounds and the best structures found so far, which needs further
investigation.

\paragraph{Use of AI.} The numerical analyses and analytical calculations, the computer code that produced the figures,
and the draft of the text were done with the large language model Claude Opus 5.5 (Anthropic),
in the course of more than 1000 prompts, under the close supervision of
the authors. The authors posed the problems, chose the approach and interpreted the results.
The authors verified the results, checked the references and edited the text, and they take full
responsibility for the content of this paper.

\section*{References}
\addcontentsline{toc}{section}{References}
\begin{small}
\begin{itemize}
\item[] Albin, N. and Cherkaev, A. (2006) Optimality conditions on fields in microstructures and
controllable differential schemes. In \emph{Inverse Problems, Multi-Scale Analysis and Effective
Medium Theory} (H. Ammari and H. Kang, eds.), Contemp. Math. \textbf{408}, Amer. Math. Soc.,
Providence, RI, 137--150.
\item[] Albin, N., Cherkaev, A. and Nesi, V. (2007) Multiphase laminates of extremal effective
conductivity in two dimensions. \emph{J. Mech. Phys. Solids} \textbf{55}, 1513--1553.
\item[] Albin, N., Conti, S. and Nesi, V. (2007) Improved bounds for composites and rigidity of
gradient fields. \emph{Proc. R. Soc. A} \textbf{463}, 2031--2048.
\item[] Alessandrini, G. and Nesi, V. (2001) Univalent $\sigma$-harmonic mappings. \emph{Arch.
Ration. Mech. Anal.} \textbf{158}, 155--171.
\item[] Avellaneda, M. (1987) Iterated homogenization, differential effective medium theory and
applications. \emph{Comm. Pure Appl. Math.} \textbf{40}, 527--554.
\item[] Ball, J. M., Currie, J. C. and Olver, P. J. (1981) Null Lagrangians, weak continuity, and
variational problems of arbitrary order. \emph{J. Funct. Anal.} \textbf{41}, 135--174.
\item[] Bauman, P., Marini, A. and Nesi, V. (2001) Univalent solutions of an elliptic system of
partial differential equations arising in homogenization. \emph{Indiana Univ. Math. J.}
\textbf{50}, 747--757.
\item[] Briane, M., Milton, G. W. and Nesi, V. (2004) Change of sign of the corrector's determinant
for homogenization in three-dimensional conductivity. \emph{Arch. Ration. Mech. Anal.}
\textbf{173}, 133--150.
\item[] Briane, M. and Nesi, V. (2004) Is it wise to keep laminating? \emph{ESAIM Control Optim.
Calc. Var.} \textbf{10}, 452--477.
\item[] Cherkaev, A. (2000) \emph{Variational Methods for Structural Optimization}. Applied
Mathematical Sciences \textbf{140}, Springer, New York.
\item[] Cherkaev, A. (2009) Bounds for effective properties of multimaterial two-dimensional
conducting composites. \emph{Mech. Mater.} \textbf{41}, 411--433.
\item[] Cherkaev, A. (2012) Optimal three-material wheel assemblage of conducting and elastic
composites. \emph{Int. J. Eng. Sci.} \textbf{59}, 27--39.
\item[] Cherkaev, A. (2024) Structure of fields in extremal 2D conducting multimaterial composites.
In \emph{Continuum Models and Discrete Systems: CMDS-14} (F. Willot, J. Dirrenberger, S. Forest,
D. Jeulin and A. V. Cherkaev, eds.), Springer Proc. Math. Stat. \textbf{457}, Springer, Cham,
75--88.
\item[] Cherkaev, A. and Gibiansky, L. V. (1996) Extremal structures of multiphase heat conducting
composites. \emph{Int. J. Solids Struct.} \textbf{33}, 2609--2623.
\item[] Cherkaev, A. and Zhang, Y. (2011) Optimal anisotropic three-phase conducting composites:
plane problem. \emph{Int. J. Solids Struct.} \textbf{48}, 2800--2813.
\item[] Gibiansky, L. V. and Sigmund, O. (2000) Multiphase composites with extremal bulk modulus.
\emph{J. Mech. Phys. Solids} \textbf{48}, 461--498.
\item[] Grabovsky, Y. (1996) Bounds and extremal microstructures for two-component composites: a
unified treatment based on the translation method. \emph{Proc. R. Soc. Lond. A} \textbf{452},
919--944.
\item[] Hashin, Z. (1962) The elastic moduli of heterogeneous materials. \emph{J. Appl. Mech.}
\textbf{29}, 143--150.
\item[] Hashin, Z. and Rosen, B. W. (1964) The elastic moduli of fiber-reinforced materials.
\emph{J. Appl. Mech.} \textbf{31}, 223--232.
\item[] Hashin, Z. and Shtrikman, S. (1962) A variational approach to the theory of the effective
magnetic permeability of multiphase materials. \emph{J. Appl. Phys.} \textbf{33}, 3125--3131.
\item[] Liu, L. P. (2011) New optimal microstructures and restrictions on the attainable
Hashin--Shtrikman bounds for multiphase composite materials. \emph{Philos. Mag. Lett.} \textbf{91},
473--482.
\item[] Liu, L. P., James, R. D. and Leo, P. H. (2007) Periodic inclusion--matrix microstructures
with constant field inclusions. \emph{Metall. Mater. Trans. A} \textbf{38}, 781--787.
\item[] Lurie, K. A. and Cherkaev, A. V. (1982) Accurate estimates of the conductivity of mixtures
formed of two materials in a given proportion (two-dimensional problem). \emph{Dokl. Akad. Nauk
SSSR} \textbf{264}, 1128--1130.
\item[] Lurie, K. A. and Cherkaev, A. V. (1984) Exact estimates of conductivity of composites formed
by two isotropically conducting media taken in prescribed proportion. \emph{Proc. R. Soc.
Edinburgh} \textbf{99A}, 71--87.
\item[] Lurie, K. A. and Cherkaev, A. V. (1985) Optimization of properties of multicomponent
isotropic composites. \emph{J. Optim. Theory Appl.} \textbf{46}, 571--580.
\item[] Lurie, K. A. and Cherkaev, A. V. (1986) Exact estimates of the conductivity of a binary
mixture of isotropic materials. \emph{Proc. R. Soc. Edinburgh} \textbf{104A}, 21--38.
\item[] Milton, G. W. (1981) Concerning bounds on the transport and mechanical properties of
multicomponent composite materials. \emph{Appl. Phys. A} \textbf{26}, 125--130.
\item[] Milton, G. W. (1986) Modelling the properties of composites by laminates. In
\emph{Homogenization and Effective Moduli of Materials and Media} (J. L. Ericksen, D. Kinderlehrer,
R. Kohn and J.-L. Lions, eds.), IMA Vol. Math. Appl. \textbf{1}, Springer, New York, 150--174.
\item[] Milton, G. W. (1990) On characterizing the set of possible effective tensors of composites:
the variational method and the translation method. \emph{Comm. Pure Appl. Math.} \textbf{43},
63--125.
\item[] Milton, G. W. (2002) \emph{The Theory of Composites}. Cambridge University Press.
\item[] Morrey, C. B. (1966) \emph{Multiple Integrals in the Calculus of Variations}. Springer.
\item[] Murat, F. and Tartar, L. (1985) Calcul des variations et homog\'en\'eisation. In \emph{Les
m\'ethodes de l'homog\'en\'eisation: th\'eorie et applications en physique}, Eyrolles, 319--369.
English translation in \emph{Topics in the Mathematical Modelling of Composite Materials}
(A. Cherkaev and R. Kohn, eds.), Progr. Nonlinear Differential Equations Appl. \textbf{31},
Birkh\"auser, Boston, 1997, 139--173.
\item[] Nesi, V. (1995) Bounds on the effective conductivity of two-dimensional composites made of
$n\ge3$ isotropic phases in prescribed volume fraction: the weighted translation method. \emph{Proc.
R. Soc. Edinburgh} \textbf{125A}, 1219--1239.
\item[] Norris, A. N. (1985) A differential scheme for the effective moduli of composites.
\emph{Mech. Mater.} \textbf{4}, 1--16.
\item[] Reshetnyak, Yu. G. (1967) On the stability of conformal mappings in multidimensional spaces.
\emph{Siberian Math. J.} \textbf{8}, 69--85.
\item[] Schulgasser, K. (1983) Sphere assemblage model for polycrystals and symmetric materials.
\emph{J. Appl. Phys.} \textbf{54}, 1380--1382.
\item[] Talbot, D. R. S., Willis, J. R. and Nesi, V. (1995) On improving the Hashin--Shtrikman
bounds for the effective properties of three-phase composite media. \emph{IMA J. Appl. Math.}
\textbf{54}, 97--107.
\item[] Tartar, L. (1979) Estimation de coefficients homog\'en\'eis\'es. In \emph{Computing Methods
in Applied Sciences and Engineering, 1977, I} (R. Glowinski and J.-L. Lions, eds.), Lecture Notes
in Math. \textbf{704}, Springer, Berlin, 364--373.
\item[] Tartar, L. (1985) Estimations fines des coefficients homog\'en\'eis\'es. In \emph{Ennio De
Giorgi Colloquium}, Pitman, 168--187.
\item[] Wiener, O. (1912) Die Theorie des Mischk\"orpers f\"ur das Feld der station\"aren
Str\"omung. Erste Abhandlung: Die Mittelwerts\"atze f\"ur Kraft, Polarisation und Energie.
\emph{Abh. Math.-Phys. Kl. K\"onigl. S\"achs. Ges. Wiss.} \textbf{32}, 509--604.
\end{itemize}
\end{small}

\end{document}